\pdfoutput=1
\documentclass[11pt]{amsart}

\usepackage[margin=1.15in]{geometry}
\usepackage[T1]{fontenc}
\usepackage{amsmath,amssymb,amsthm}
\usepackage{graphicx}
\usepackage{booktabs}
\usepackage{enumitem}
\usepackage{tikz}
\usepackage{newunicodechar}
\usepackage{needspace}
\usepackage{etoolbox}
\newunicodechar{∀}{\ensuremath{\forall}}
\newunicodechar{∃}{\ensuremath{\exists}}
\newunicodechar{∈}{\ensuremath{\in}}
\newunicodechar{ℝ}{\ensuremath{\mathbb{R}}}
\newunicodechar{×}{\ensuremath{\times}}
\newunicodechar{ˢ}{\ensuremath{^{s}}}
\newunicodechar{∧}{\ensuremath{\wedge}}
\newunicodechar{₀}{\ensuremath{_0}}
\usepackage[colorlinks=true,linkcolor=blue!50!black,citecolor=blue!50!black,urlcolor=blue!50!black]{hyperref}

\newtheorem*{mainthm}{Main Theorem}
\newtheorem{thm}{Theorem}[section]
\newtheorem{lemma}[thm]{Lemma}
\newtheorem{proposition}[thm]{Proposition}
\theoremstyle{definition}
\newtheorem{definition}[thm]{Definition}
\theoremstyle{remark}
\newtheorem{remark}[thm]{Remark}

\newcommand{\R}{\mathbb{R}}
\newcommand{\magn}{\operatorname{mag}}
\newcommand{\mig}{\operatorname{mig}}
\newcommand{\ival}[1]{[#1]}
\newcommand{\eu}{e_u}
\newcommand{\es}{e_s}
\newcommand{\Wu}{W^u}
\newcommand{\Ws}{W^s}
\newcommand{\loc}{\mathrm{loc}}
\newcommand{\cs}{c_*}
\newcommand{\cmid}{c_{\mathrm{mid}}}
\newcommand{\Kr}{\mathcal{K}}
\newcommand{\Dh}{\mathcal{D}}

\begin{document}

\title[SRB measures for some R\"ossler systems]{Existence of SRB measures for some R\"ossler systems}

\author{Theophile Caby}
\address{Laboratoire d'Oc\'eanographie Physique et Spatiale, Univ Brest, CNRS, Ifremer, IRD, Brest, France}
\email{theophile.caby@univ-brest.fr}
\makeatletter
\apptocmd{\@setauthors}{%
  \par\vspace{0.6\baselineskip}%
  {\centering\normalfont\small
   Laboratoire d'Oc\'eanographie Physique et Spatiale, Univ Brest, CNRS, Ifremer, IRD, Brest, France\\
   \texttt{theophile.caby@univ-brest.fr}\\[0.3\baselineskip]
   Code and logs: \url{https://github.com/cabytheo-cloud/Rossler_SRB}\par}%
}{}{\PackageWarning{paper}{could not patch author block}}
\makeatother

\date{\today}


\begin{abstract}
We prove that for the R\"ossler system $\dot x=-(y+z)$, $\dot y=x+ay$, $\dot z=b+z(x-c)$ with $a=b=0.2$, the
set of parameters $c$ within $6\cdot10^{-5}$ of the classical value $c=5.7$ for which the flow has a
H\'enon-like strange attractor supporting a Sinai--Ruelle--Bowen (SRB) measure has positive Lebesgue
measure. The key step, proved with computer assistance, is that at a parameter localised in an explicit interval within
$5.5\cdot10^{-5}$ of $5.7$ the period-one saddle of the Poincar\'e map has a quadratic homoclinic tangency
that unfolds generically; the conclusion then follows from the theorems of Mora--Viana, Benedicks--Young and
Benedicks--Viana. The same method gives two further such tangencies near $5.7$ and one at the band-merging
parameter $c\approx4.507$. The proof reduces the generic-tangency conditions to finitely many inequalities,
through explicit $C^2$ graph-transform bounds for the local invariant manifolds of a Poincar\'e map whose
area contraction is of order $10^{-11}$--$10^{-14}$, and a scalar tangency function. The inequalities are
checked with interval arithmetic in CAPD and re-checked with an independent implementation based on the Kv
library; the analytic tangency criterion is formalised in Lean~4. Whether an SRB measure exists at $c=5.7$
itself remains open.
\end{abstract}

\maketitle

\section{Introduction}

The R\"ossler system~\cite{Rossler76}
\begin{equation}\label{eq:rossler}
\dot x=-(y+z),\qquad \dot y=x+ay,\qquad \dot z=b+z(x-c)
\end{equation}
is one of the simplest three-dimensional flows with an apparently chaotic attractor. At the classical
parameters $(a,b,c)=(0.2,0.2,5.7)$, numerical simulations show a strange attractor whose first return
map to a transversal section is, to all visible precision, a unimodal interval map. A natural
mathematical question is whether this attractor carries a \emph{Sinai--Ruelle--Bowen (SRB) measure}: an
invariant probability with a positive Lyapunov exponent whose conditional measures along unstable
manifolds are absolutely continuous. Such a measure describes the statistics of a set of initial
conditions of positive volume (see~\cite{Young02} for a survey).

For uniformly hyperbolic and singular-hyperbolic attractors, SRB measures are well understood. In
particular, Tucker's computer-assisted proof~\cite{Tucker02} showed that the Lorenz equations have a
singular-hyperbolic attractor, and such attractors support an SRB measure~\cite{APPV09}. The decisive
feature there is that singular hyperbolicity is an \emph{open} condition, defined by invariant cone
fields, which can be certified by a finite computation.

The R\"ossler attractor is of a different nature. Its return map folds, so no invariant cone field can
exist near the fold; the attractor is at best non-uniformly hyperbolic, of H\'enon type. For such
systems, horseshoes of positive entropy are robust, but the existence of SRB measures is a delicate,
parameter-dependent property. For the H\'enon family it has been proved only for positive-measure,
Cantor-like sets of parameters~\cite{BC91,BY93,MV93,BV01}, and the question remains open at the
classical H\'enon parameters.

Moreover, spiral attractors of R\"ossler type are \emph{quasi-attractors} in the sense of
Afraimovich--Shilnikov: stable periodic orbits of long period are present, or appear under arbitrarily
small perturbations~\cite{GGKK18}. Galias~\cite{Gal26} recently fixed $c=5.7$, varied $a=b$ jointly, and
proved with interval arithmetic that the system has a stable periodic orbit (of period $97$ for the return
map) for parameters in a periodic window of width about $4\cdot10^{-38}$ at distance less than
$2\cdot10^{-22}$ from the classical value $a=b=0.2$. His results indicate that periodic windows are dense
near the classical parameters, and the widths of the windows he found provide, in his words, a strong
argument that the classical attractor itself is not periodic. In any case, ``having an SRB measure'' is
not an open condition in parameter space, and it cannot be certified at a given parameter value by
verifying an open condition, as was done for the Lorenz attractor.

Rigorous results on the R\"ossler system so far concern topological dynamics and periodic orbits:
\begin{itemize}[leftmargin=*]
\item symbolic dynamics and positive topological entropy at the classical parameters~\cite{Zgl97};
\item a trapping region for the return map~\cite{CAPD21};
\item period doubling~\cite{WZ09pd};
\item periodic orbits of all periods at $c=5.25$, and of all periods except $3$ at $c=4.7$~\cite{GZ21};
\item for $b=0.2$, $c=15$ and $a\in[0.12,0.3659]$, a trapping region of the Poincar\'e map containing a
connected attractor of positive topological entropy, whose entropy grows through a sequence of
bifurcations of periodic orbits~\cite{WSB25};
\item periodic windows extremely close to the classical parameters~\cite{Gal26};
\item a reduction of the non-wandering set to the quadratic family~\cite{Igra23}.
\end{itemize}
The ``chaos'' established in~\cite{Zgl97,WSB25} is topological: these works prove the existence of invariant
sets carrying symbolic dynamics, hence positive topological entropy. Such sets are typically of Lebesgue measure
zero and may coexist with a periodic sink, so they do not decide whether the attractor observed in simulations is
chaotic in the metric sense (positive Lyapunov exponent for typical orbits, or an SRB measure); that question is
open at the classical parameters, see~\cite{Gal26}. To our knowledge, no SRB measure had been established for any
parameter of the R\"ossler family.

\subsection*{Main result}

\begin{mainthm}
There is a parameter $\cs$ with $|\cs-5.7|<5.47\cdot10^{-5}$ such that, for every $\varepsilon>0$, the set of
$c\in(\cs-\varepsilon,\cs+\varepsilon)$ for which the R\"ossler flow~\eqref{eq:rossler} with $a=b=0.2$ has a
H\'enon-like strange attractor supporting an ergodic SRB measure, whose basin has positive Lebesgue measure
in $\R^3$, has positive Lebesgue measure. In particular, the set of such parameters $c$ in
$(5.7-6\cdot10^{-5},\,5.7+6\cdot10^{-5})$ has positive Lebesgue measure.
\end{mainthm}

From now on $a=b=0.2$, and $P_c$ denotes the Poincar\'e map of~\eqref{eq:rossler} on the section
$\Sigma=\{x=0,\ \dot x>0\}$, written in the coordinates $(y,z)$ (see Section~\ref{sec:setting} for the
precise definition). The parameter $\cs$ is explicitly localised: it lies in the interval $J_3\subset(5.70005464309,\,5.70005464312)$
of~\eqref{eq:Ji}. The proof has two parts of different nature.
\begin{itemize}[leftmargin=*]
\item \emph{The new ingredient}, proved with computer assistance (Proposition~\ref{prop:tangency}): at
$c=\cs$ the period-one saddle $p_c$ of $P_c$ has a quadratic homoclinic tangency that unfolds generically
as $c$ varies. In the (almost one-dimensional) return map, this corresponds to the critical point landing
on $p$ after $16$ iterates (Figure~\ref{fig:returnmap}).
\item \emph{The passage from such a tangency to SRB measures} (Proposition~\ref{prop:srb}), which relies on
the theorem of Mora and Viana~\cite{MV93} on the abundance of strange attractors in generic unfoldings of
homoclinic tangencies, on the construction of SRB measures for these attractors and the solution of the
basin problem~\cite{BY93,BV01,Viana98}, and on a standard suspension argument. One hypothesis
of~\cite{MV93}, the smooth linearisability of the saddle, reduces to a finite-order non-resonance condition on its
multipliers, which we verify with interval arithmetic (Remark~\ref{rem:MVhyp}).
\end{itemize}
The same argument gives positive-measure sets of SRB parameters accumulating at a further tangency
parameter, at distance less than $5.81\cdot10^{-4}$ from $5.7$ (Proposition~\ref{prop:tangency}), and at the
parameter $c\approx4.507$ where the two bands of the chaotic attractor of the period-doubling cascade merge
(Proposition~\ref{prop:band}). A third tangency near $5.7$, at distance less than $4.21\cdot10^{-4}$, was verified
with CAPD only, and there the non-resonance condition has been checked only up to order $28$.

These results say nothing about the classical value $c=5.7$ itself, nor about the tangency parameter $\cs$.
Nor do they identify the attractors they
produce with the large attractor seen in simulations: they are the small H\'enon-like attractors created
near the orbit of tangency. We return to these points in Section~\ref{sec:discussion}.

\subsection*{Method}
A quadratic, generically unfolding tangency is characterised by four conditions on a scalar function
$g(s,c)$, which measures the signed distance from a parametrised piece of the unstable manifold to the
local stable manifold:
\[
g=0,\qquad \partial_s g=0,\qquad \partial_s^2 g\neq 0,\qquad \partial_c g\neq 0 .
\]
To make these conditions checkable, we need the local invariant manifolds of $p_c$ with
\emph{second-order} control, including their dependence on $c$. We obtain this from explicit $C^2$
graph-transform estimates (Lemmas~\ref{lem:Wu} and~\ref{lem:Ws}), which reduce everything to interval
inequalities on the first and second derivatives of $P_c$ over small boxes. The remaining quantities
are derivatives of $P_c^N$ along the unstable manifold, which we enclose with the $C^2$-Lohner
algorithm of the CAPD library~\cite{CAPD21,WZ11}. Each verification takes a few seconds on one core.

The approach is close in spirit to the computer-assisted proof of a homoclinic tangency for the H\'enon
map and the forced damped pendulum by Wilczak and Zgliczy\'nski~\cite{WZ09}, which uses covering
relations and cone conditions in the projectivised tangent bundle; see also~\cite{AM06}. Coomes, Ko\c
cak and Palmer~\cite{CKP26} recently gave a computable criterion for generic quadratic tangencies of a
fixed point of a planar diffeomorphism depending on two parameters. It is based on Kleinkauf's
characterisation of quadratic, generically unfolding tangencies through bounded solutions of the
variational equations along the homoclinic orbit, and reduces the proof to the invertibility of a finite
matrix built from a pseudo-orbit, via the implicit function theorem; it is illustrated on explicit
planar maps (H\'enon, Duffing, Bogdanov, Perennial Grass). Our setting and method are different. The
map is a Poincar\'e map of a three-dimensional flow, known only through validated integration, and it is
extremely dissipative: $|\lambda_s|\approx10^{-11}$--$10^{-14}$, so the stable manifold is very well
determined but the unstable fold is correspondingly thin. What our method provides is:
\begin{itemize}[leftmargin=*]
\item explicit $C^2$ bounds on both local manifolds, together with their parameter dependence, from
graph transforms whose contraction factors are about $2\cdot10^{-4}$ for $\Wu$ and $0.42$ for $\Ws$;
the extreme dissipation is an advantage here rather than a difficulty;
\item correlated (affine) enclosures of $P_c^N$ for $N$ up to $26$;
\item the generic-unfolding condition $\partial_cg\neq0$ in a one-parameter family;
\item the passage, via Mora--Viana, to SRB measures for the flow.
\end{itemize}

\subsection*{Independent checks}
To reduce the trust placed in a single code base, we did two further checks.
\begin{itemize}[leftmargin=*]
\item \emph{Second implementation.} The hypotheses for $J_1$, $J_3$ and $J_0$, and the numerical
statements of Propositions~\ref{prop:tangency} and~\ref{prop:band} except the bound on the distance from the tangency
point to $p$, were re-verified with a separately written program
built on a different validated-numerics library, Kv~\cite{Kv} (Section~\ref{sec:kv}). It computes the
Poincar\'e map in a different way: as a fixed-angle map after reparametrising the flow by the polar
angle.
\item \emph{Formal proof.} The analytic tangency criterion (Lemma~\ref{lem:criterion}) and the Taylor
bound (Remark~\ref{rem:taylor}) were formalised and machine-checked in Lean~4 with
Mathlib~\cite{Lean4,mathlib} (Section~\ref{sec:lean}).
\end{itemize}

\subsection*{Organisation}
Section~\ref{sec:setting} fixes definitions. Section~\ref{sec:B} shows that a generically unfolding
tangency gives SRB measures (Proposition~\ref{prop:srb}), listing the hypotheses of the theorems used.
Section~\ref{sec:reduction} reduces the existence of such a tangency to finitely many inequalities, which
are verified in Section~\ref{sec:cap} near $c=5.7$, where the Main Theorem is proved, and in
Section~\ref{sec:band} at the band-merging parameter. Section~\ref{sec:indep} describes the independent
implementation, the Lean formalisation and the code audit. Section~\ref{sec:discussion} discusses the
tangency parameter $\cs$ and the classical parameter.

\section{Setting and definitions}\label{sec:setting}

\subsection{The Poincar\'e map}
Let $X_c$ denote the vector field~\eqref{eq:rossler} with $a=b=0.2$, and let $\varphi^t_c$ be its flow.
The plane $\{x=0\}$ is crossed transversally in the positive direction exactly on
$\Sigma=\{(0,y,z):\ y+z<0\}$. For $q\in\Sigma$ whose forward orbit returns to $\Sigma$, let
$\tau_c(q)>0$ be the first return time and $P_c(q)=\varphi_c^{\tau_c(q)}(q)$. On the open set where it is
defined, $P_c$ is an injective local diffeomorphism, real analytic jointly in $(q,c)$, because $X_c$ is
polynomial and the crossing is transversal. We use $(y,z)$ as coordinates on $\Sigma$.

It is convenient to treat $c$ as an additional state variable. The extended field
$\hat X(x,y,z,c)=(X_c(x,y,z),0)$ has a Poincar\'e map $\hat P(y,z,c)=(P_c(y,z),c)$ on $\{x=0\}$, whose first
and second derivatives contain all the $c$-derivatives we need. These are derivatives of the Poincar\'e
map itself, including the correction for the varying return time: with $e_1=(1,0,0)$,
\[
DP_c(q)=\Bigl(I-\frac{X_c(P_c(q))\,e_1^T}{\dot x(P_c(q))}\Bigr)D\varphi_c^{\tau_c(q)}(q)
\quad\text{restricted to } T_q\Sigma,
\]
and similarly at second order and in the $c$-direction; Section~\ref{sec:audit} explains how they are
computed.

\subsection{Orientation}
Since $\Sigma$ is crossed in one direction only, $P_c$ preserves orientation. Explicitly, for $q$ in its
domain,
\[
\det DP_c(q)=\exp\Bigl(\int_0^{\tau_c(q)}\operatorname{div}X_c(\varphi^t_cq)\,dt\Bigr)\,
\frac{\dot x(q)}{\dot x(P_c(q))}>0 ,
\]
because $\dot x=-(y+z)>0$ on $\Sigma$. (Compare the volumes spanned by $X_c,e_y,e_z$ at $q$ and their images
under $D\varphi_c^{\tau_c(q)}$, and use Liouville's formula.) At a fixed point the last factor equals $1$.
In particular, if $p$ is a saddle fixed point with $\lambda_u<-1$, then
$\lambda_s=\det DP_c(p)/\lambda_u<0$.

\subsection{Homoclinic tangencies}
Let $(F_c)_{c\in I}$ be a $C^2$ family of injective local diffeomorphisms of an open subset of $\R^2$,
and let $p_c$ be a hyperbolic saddle fixed point of $F_c$ depending $C^2$ on $c$.

\begin{definition}\label{def:tangency}
We say that $F_{\cs}$ has a \emph{quadratic homoclinic tangency associated to $p_{\cs}$, unfolding
generically at $\cs$}, if there exist
\begin{itemize}[leftmargin=*]
\item an open set $\mathcal U\subset\R^2$ and an interval $I'\ni\cs$;
\item a $C^2$ function $\ell\colon\mathcal U\times I'\to\R$ with $\nabla_q\ell\neq0$, such that for each
$c$ the zero set $\{\ell(\cdot,c)=0\}$ is an arc of $\Ws(p_c)$;
\item a $C^2$ map $\Gamma\colon(s_0-\eta,s_0+\eta)\times I'\to\mathcal U$ with $\partial_s\Gamma\neq0$ and
$\Gamma(s,c)\in\Wu(p_c)$;
\end{itemize}
such that the function $g(s,c)=\ell(\Gamma(s,c),c)$ satisfies
\[
g(s_0,\cs)=0,\qquad \partial_s g(s_0,\cs)=0,\qquad \partial_s^2 g(s_0,\cs)\neq0,\qquad
\partial_c g(s_0,\cs)\neq0 .
\]
\end{definition}

\begin{remark}\label{rem:normalform}
The conditions are independent of the choice of the defining function $\ell$ and of the regular
parametrisation $\Gamma$. Straightening $\Ws(p_c)$ to $\{x_2=0\}$, the unstable arc becomes, near the
tangency point, a graph $x_2=\alpha x_1^2+\beta(c-\cs)+\dots$, with $\alpha\neq0$ equivalent to
$\partial_s^2g\neq0$ and $\beta\neq0$ to $\partial_c g\neq0$. These are the two conditions of a generic
unfolding of a homoclinic tangency in the sense of Palis--Takens~\cite{PT93} and of
Mora--Viana~\cite[p.~6]{MV93}: in the normal form of~\cite{MV93},
$f^N_\mu(1+\xi,\eta)=(\alpha\xi^2+\beta\eta+v\mu+\dots,\,1+\dots)$ with $\alpha\neq0$ and $v\neq0$. See
also~\cite[Def.~4]{WZ09}.
\end{remark}

\subsection{SRB measures and strange attractors}
Let $F$ be a $C^2$ diffeomorphism (or flow) and $\mu$ an ergodic invariant Borel probability measure
with compact support. We call $\mu$ an \emph{SRB measure} if it has a positive Lyapunov exponent
$\mu$-almost everywhere (for flows, one transverse to the flow direction), and its conditional measures
on local unstable manifolds are absolutely continuous with respect to the induced Riemannian measure.
The \emph{basin} of $\mu$ is the set of points $q$ whose empirical measures
$\frac1n\sum_{k<n}\delta_{F^kq}$ (respectively $\frac1T\int_0^T\delta_{\varphi^tq}\,dt$) converge
weakly$^*$ to $\mu$.

Following~\cite[p.~2]{MV93}, an \emph{attractor} of a surface map is a compact invariant set $\Lambda$ with a
dense orbit whose stable set $\Ws(\Lambda)$ (its topological basin) has non-empty interior; it is
\emph{strange} if it contains a dense orbit with a positive Lyapunov exponent. We call \emph{H\'enon-like}
the strange attractors produced by the renormalisation of~\cite{MV93} near a homoclinic tangency
(Theorem~\ref{thm:MV}); they are not uniformly hyperbolic~\cite[p.~2]{MV93}. For the R\"ossler flow, a
H\'enon-like strange attractor is the flow saturation $\{\varphi^t_cq: q\in\Lambda,\ t\in\R\}$ of a
H\'enon-like strange attractor $\Lambda\subset\Sigma$ of $P_c$. By the flow-box theorem, its basin also has
non-empty interior.

\section{From homoclinic tangencies to SRB measures}\label{sec:B}

We use three results on generic unfoldings of homoclinic tangencies of surface maps. We state them in
the form in which we use them, with precise references.

\begin{thm}[Mora--Viana~{\cite[Theorem~A and \S2]{MV93}}]\label{thm:MV}
Let $(f_\mu)_\mu$ be a $C^\infty$ one-parameter family of surface diffeomorphisms such that $f_0$ has a
homoclinic tangency associated to a hyperbolic fixed point $p_0$, with the following properties
(\cite[pp.~5--6]{MV93}):
\begin{enumerate}[label=\textup{(MV\arabic*)},leftmargin=*]
\item\label{MV1} the saddle is dissipative: $|\det Df_\mu(p_\mu)|<1$;
\item\label{MV2} its eigenvalues satisfy $0<\lambda<1<\sigma$;
\item\label{MV3} the eigenvalues satisfy $\lambda^\ell\sigma^m\neq1$ for all integers $\ell,m\ge1$ with $\ell+m\le r-2$,
for a fixed, sufficiently large $r$ ($r=140$ suffices); this gives $\mu$-dependent linearising coordinates for
$f_\mu$ of class $C^{[(r-1)/2]}$ near $p_\mu$;
\item\label{MV4} the tangency is quadratic and unfolds generically with $\mu$ (Remark~\ref{rem:normalform}).
\end{enumerate}
Then, for parameter intervals $I_n$ converging to $\mu=0$ and neighbourhoods $U_n$ of a point of the orbit of
tangency, suitable rescalings of $f_\mu^{\,n+N}|U_n$, $\mu\in I_n$, converge in the $C^k$ topology, for some $k$ growing with $r$, to the
quadratic family $(a,x,y)\mapsto(1-ax^2,0)$ and form a H\'enon-like family in the sense
of~\cite[Theorem~2.1]{MV93}; and for $n$ large there is a set $E_n\subset I_n$ of positive Lebesgue measure such
that, for $\mu\in E_n$, the map $f_\mu^{\,n+N}|U_n$ has a strange attractor~\cite[Theorem~B]{MV93}. In
particular, for every $\varepsilon>0$ the set of $\mu\in(-\varepsilon,\varepsilon)$ for which $f_\mu$ has a
strange attractor near the orbit of tangency has positive Lebesgue measure.
\end{thm}

\begin{remark}[on the hypotheses]\label{rem:MVhyp}
We comment on each hypothesis, for the family $f_\mu=P^2_{\cs+\mu}$ used in the proof of
Proposition~\ref{prop:srb}.
\begin{itemize}[leftmargin=*]
\item \emph{Regularity.} $(P_c)$ is real analytic in $(q,c)$.
\item \ref{MV1}, \ref{MV2}: $\det DP^2_c(p_c)=(\lambda_u\lambda_s)^2<1$, and the eigenvalues of $DP^2_c(p_c)$ are
$\lambda_s^2\in(0,1)$ and $\lambda_u^2>1$. (Mora and Viana work with positive eigenvalues, which can always
be achieved by passing to the second iterate, as we do for our flip saddle.)
\item \ref{MV4} is exactly what Definition~\ref{def:tangency} encodes; $P^2$ has the same invariant manifolds,
tangency and function $g$ as $P$.
\item \ref{MV3} is the form in which Mora and Viana state their linearisation hypothesis~\cite[\S2]{MV93}; we
follow the detailed account of Hiraide~\cite[Theorem~2 and \S1]{Hir92}, which gives the condition as above, notes
that $r=140$ suffices, and obtains the $C^{[(r-1)/2]}$ linearising coordinates, depending smoothly on $\mu$, from
Gomozov's theorem (see~\cite{Bel78}). For a saddle with $0<\lambda<1<\sigma$ a relation $\lambda^\ell\sigma^m=1$ means
that $\rho:=\log(1/\lambda)/\log\sigma$ equals $m/\ell$; passing from $P$ to $P^2$ does not change $\rho$. The Kv
programs of Section~\ref{sec:kv} enclose $\rho$ for all $c$ in $J_1$, $J_3$ and $J_0$ and certify that $\rho\ne p/q$
for all integers $p,q\ge1$ with $p+q+1\le10^3$, which contains the range $\ell+m\le138$ required for $r=140$.
(A relation $\rho=p/q$ is equivalent to the resonances $\lambda=\lambda^{q+1}\sigma^p$ and $\sigma=\lambda^q\sigma^{p+1}$, of
order $p+q+1$; we accordingly say that there is no resonance of order at most $10^3$.) In
fact, for $J_3$ they give $\rho\in[36.4692539857,\,36.4692539884]$, and no relation $\rho=p/q$ with
$p+q+1<1\,229\,030$ holds; for $J_1$ and $J_0$ the corresponding bounds are $1\,200$ (the Kv enclosure of
$\det DP_c(p_c)$ is much wider on $J_1$) and $45\,389$. (Renormalisation schemes that avoid linearisation altogether
were developed by Romero~\cite{Romero95} and by Gonchenko, Shilnikov and Turaev~\cite{GST08}; we do not use them.)
\item Barrientos and Rojas report~\cite[footnote~2]{BR21} that the proof of~\cite{MV93} needs an additional
assumption on the distortion of the determinant of the family, which seems unnecessary in view of Berger's
alternative approach~\cite{Ber19}. Under~\ref{MV3} this assumption holds automatically. In linearising coordinates the return maps near the tangency are $T_1\circ L_\mu^n$ with $L_\mu$
linear, so their Jacobian is $(\lambda\sigma)^n\det DT_1$, where $\det DT_1$ is a fixed smooth function that does
not vanish; after rescaling it is evaluated on a neighbourhood of the tangency point that shrinks
exponentially, so the relative variation of the Jacobian tends to $0$. Berger's version of the Palis--Takens
renormalisation theorem~\cite[Theorem~1.4 and Remark~1.2]{BerNF}, whose proof uses only $C^r$ linearising
coordinates, makes this precise: for every $\delta>0$ and $n$ large, the renormalised family has the form
$(x,y)\mapsto(x^2+a+y,\,-bx)+(A_a,\,b\,B_a)(x,y)$ with $A$ and $B$ $\delta$-small in $C^r$, and its determinant lies
in $[(1-3\delta)b,(1+3\delta)b]$.
\item Wang and Young~\cite[Theorem~9]{WY01} state the passage from a generic unfolding to SRB measures with
explicit hypotheses: a $C^\infty$ family unfolding a non-degenerate tangency of a hyperbolic fixed point with
$0<\lambda<1<\sigma$ and $\lambda\sigma<1$, whose multipliers lie in the open and dense set of pairs to which
Sternberg's theorem applies. These are \ref{MV1}, \ref{MV2}, \ref{MV4} and a finite-order non-resonance
condition of the type discussed above. Their conclusions (\cite[Theorems~1--7]{WY01}) contain those of
Theorem~\ref{thm:BV} below for the SRB measure and its basin, and add that there is at most one ergodic SRB
measure (the limiting interval map has a single critical point), exponential decay of correlations and a
central limit theorem.
\end{itemize}
\end{remark}

\begin{thm}[Benedicks--Young~\cite{BY93}, Benedicks--Viana~\cite{BV01}; see~{\cite[\S4, Theorems~3 and~4]{Viana98}} and~{\cite[Theorems~3 and~9]{WY01}}]\label{thm:BV}
In the situation of Theorem~\ref{thm:MV}, after replacing each $E_n$ by a subset of positive Lebesgue
measure, the strange attractor $\Lambda$ of $f_\mu^{\,n+N}|U_n$, $\mu\in E_n$, given by Theorem~\ref{thm:MV} supports a unique SRB measure, and
Lebesgue-almost every point of the topological basin of $\Lambda$ belongs to the basin of this measure.
\end{thm}

Benedicks and Young~\cite{BY93} constructed SRB measures for the H\'enon attractors of Benedicks and
Carleson~\cite{BC91}. Benedicks and Viana~\cite{BV01} proved, for a large class of non-uniformly hyperbolic
attractors of dissipative diffeomorphisms, that there are no ``holes'' in the basin: Lebesgue-almost every
point of the topological basin is generic for the SRB measure. Their analysis is carried out in the
setting of the H\'enon-like families of Mora and Viana~\cite[\S3]{Ben02}, and Viana~\cite[\S4]{Viana98}
states both results for H\'enon-like attractors, including those arising in generic unfoldings of
homoclinic tangencies, under the strong area dissipativeness that holds for the renormalised maps of
Theorem~\ref{thm:MV}. Mora and Viana had anticipated this extension~\cite[p.~3]{MV93}.

Both theorems are local: all the maps and iterates they use are defined on a neighbourhood of the closure
of the orbit of tangency, on which $P_c$ is an analytic diffeomorphism onto its image, depending
analytically on $c$. They therefore apply to the family $(P_c)$.

\begin{lemma}[suspension]\label{lem:suspension}
Let $\Lambda\subset\Sigma$ be a compact $P_c$-invariant attractor in the domain of $P_c$, and let $\nu$ be
an ergodic SRB measure of $P_c$ supported on $\Lambda$ whose basin meets the topological basin of
$\Lambda$ in a set of positive area. Then
\[
\mu=\frac{1}{\int\tau_c\,d\nu}\int_\Lambda\int_0^{\tau_c(q)}\delta_{\varphi_c^t q}\,dt\,d\nu(q)
\]
is an ergodic SRB measure for the flow $\varphi_c^t$, supported on the flow saturation of $\Lambda$, and
its basin has positive Lebesgue measure in $\R^3$.
\end{lemma}

\begin{proof}
This is standard; we indicate the main points. Let $V$ be a compact neighbourhood of $\Lambda$ in
$\Sigma$ on which $\tau_c$ is continuous and bounded away from $0$ and $\infty$. Invariance and
ergodicity are classical for suspension flows. The exponents of $\mu$ transverse to the flow are those of
$\nu$ divided by $\int\tau_c\,d\nu$; in particular one is positive. The strong unstable manifold of
$\varphi_c$ through $q\in\Lambda$ projects diffeomorphically, along flow lines, onto $\Wu(q;P_c)$, so
absolute continuity of conditional measures is preserved. For the basin, let $q\in B(\nu)$ lie in the
topological basin of $\Lambda$, so that $P_c^kq\in V$ for all large $k$. For continuous $\psi$, the flow
averages of $\psi$ along $q$ are, up to an error bounded by $\sup|\psi|\sup_V\tau_c/T$, ratios of Birkhoff
sums of the continuous functions $\Psi(q)=\int_0^{\tau_c(q)}\psi(\varphi^t_cq)\,dt$ and $\tau_c$, which
converge to $\int\Psi\,d\nu\big/\!\int\tau_c\,d\nu=\int\psi\,d\mu$. Hence the basin of $\mu$ contains the
flow saturation of this set, which has positive volume because $(q,t)\mapsto\varphi^t_cq$ is a local
diffeomorphism from $\Sigma\times(0,t_0)$ to $\R^3$ for small $t_0>0$.
\end{proof}

\begin{proposition}\label{prop:srb}
Let $\cs$ be a parameter at which $P_{\cs}$ has a saddle fixed point $p_{\cs}$ with
$|\det DP_{\cs}(p_{\cs})|<1$, whose multipliers satisfy the non-resonance condition~\ref{MV3}
(Remark~\ref{rem:MVhyp}), and whose stable and unstable manifolds have a quadratic homoclinic tangency
unfolding generically at $\cs$ (Definition~\ref{def:tangency}). Then, for every $\varepsilon>0$, the set of
$c\in(\cs-\varepsilon,\cs+\varepsilon)$ for which the R\"ossler flow~\eqref{eq:rossler} with $a=b=0.2$ has a
H\'enon-like strange attractor supporting an ergodic SRB measure, whose basin has positive Lebesgue measure
in $\R^3$, has positive Lebesgue measure.
\end{proposition}

\begin{proof}
By Section~\ref{sec:setting}, $0<\det DP_{\cs}(p_{\cs})<1$, so the two multipliers of the saddle have the
same sign, and $f_\mu=P_{\cs+\mu}^2$ has positive multipliers $\lambda_s^2\in(0,1)$ and $\lambda_u^2>1$. We
apply Theorems~\ref{thm:MV} and~\ref{thm:BV} to $(f_\mu)$, whose hypotheses are discussed in
Remark~\ref{rem:MVhyp}.

Fix $\varepsilon>0$. We obtain a set $\mathcal E\subset(\cs-\varepsilon,\cs+\varepsilon)$ of positive measure
such that, for $c\in\mathcal E$, an iterate $F=P_c^{m}$ (with $m=2(n+N)$) restricted to a small open set $U$
has a strange attractor $\Lambda$ supporting an SRB measure $\nu_0$ as in Theorem~\ref{thm:BV}. Let
$\Lambda'=\bigcup_{j<m}P_c^j(\Lambda)$. It is compact and $P_c$-invariant; the $P_c$-orbit of a point whose
$F$-orbit is dense in $\Lambda$ is dense in $\Lambda'$, and its Lyapunov exponent for $P_c$ is that for $F$
divided by $m$, hence positive; and the stable set of $\Lambda'$ contains that of $\Lambda$, which has non-empty
interior. So $\Lambda'$ is a strange attractor of $P_c$. The measure $\nu=\frac1m\sum_{j<m}(P_c^j)_*\nu_0$ is
an ergodic SRB measure of $P_c$ supported on $\Lambda'$: push-forwards by diffeomorphisms preserve absolute
continuity along unstable manifolds. Its basin contains the basin of $\nu_0$ for $F$, which by
Theorem~\ref{thm:BV} has full Lebesgue measure in the topological basin of $\Lambda$, a set with non-empty
interior. Lemma~\ref{lem:suspension} concludes the proof.
\end{proof}

\section{Reduction to finitely many inequalities}\label{sec:reduction}

\subsection{Interval notation and local coordinates}\label{sec:notation}
For a compact interval $\ival{x}$ we write $\magn\ival{x}=\max\{|x|:x\in\ival{x}\}$ and
$\mig\ival{x}=\min\{|x|:x\in\ival{x}\}$. Arithmetic on intervals is the usual set-valued arithmetic with
outward rounding; an inclusion such as $\ival{x}\subset[-L,L]$ is checked on the computed endpoints.

Throughout this section, $J=[c_L,c_R]$ is one of the parameter intervals and $E=[\eu\ \es]$ is a fixed
invertible $2\times2$ matrix, an approximate eigenbasis of $DP_c(p_c)$. The point $p_c$ is the fixed point
of $P_c$ given by Lemma~\ref{lem:krawczyk} below, which also provides a box $\Kr_J\ni p_c$ for all
$c\in J$. The \emph{local map}
\[
f_c(u,v)=E^{-1}\bigl(P_c(p_c+E(u,v)^T)-p_c\bigr)
\]
satisfies $f_c(0,0)=0$ for every $c\in J$. Given a box $B\subset\R^2$ in local coordinates, we say that
intervals $\ival{a_{ij}}$, $\ival{h^i_{jk}}$, $\ival{\partial_cf^i}$ \emph{enclose the derivatives of $f$
on $B\times J$} if, for all $(u,v)\in B$ and $c\in J$,
\[
\partial_j f^i_c(u,v)\in\ival{a_{ij}},\qquad
\partial_j\partial_k f^i_c(u,v)\in\ival{h^i_{jk}},\qquad
\partial_c\bigl(f^i_c(u,v)\bigr)\in\ival{\partial_cf^i},
\]
where $\partial_1=\partial_u$ and $\partial_2=\partial_v$. Differentiating the definition of $f_c$ gives
\begin{equation}\label{eq:fc}
\partial_c f_c(u,v)=E^{-1}\bigl(\partial_cP_c(q)+(DP_c(q)-I)\,p_c'\bigr),\qquad q=p_c+E(u,v)^T .
\end{equation}

\emph{Where the enclosures are computed.} The fixed point $p_c$ is only known to lie in the box $\Kr_J$.
All enclosures are therefore computed on the interval hull of $(\Kr_J+E\cdot B)\times J$, and \emph{not}
around an approximate fixed point: we enclose $DP_c$, $D^2P_c$ and $\partial_cP_c$ on this set, and
$p_c'=dp_c/dc$ on $\Kr_J\times J$. Since this set contains $\{(p_c+E(u,v)^T,c):(u,v)\in B,\ c\in J\}$, the
resulting intervals enclose the derivatives of $f$ on $B\times J$. The matrix $E^{-1}$ is replaced by a
rigorous interval enclosure computed by the adjugate formula.

\subsection{The fixed point}

\begin{lemma}[Krawczyk~\cite{Krawczyk69}]\label{lem:krawczyk}
Let $Q_\Kr\subset\R^2$ be a box, $\tilde q\in Q_\Kr$, and $R$ a $2\times2$ matrix. Suppose that, for an
interval $\ival{c}$,
\[
\Kr=\tilde q-R\bigl(\ival{P_{\ival c}(\tilde q)}-\tilde q\bigr)+\bigl(I-R(\ival{DP_{\ival c}(Q_\Kr)}-I)\bigr)(Q_\Kr-\tilde q)
\]
is contained in the interior of $Q_\Kr$, where $\ival{P_{\ival c}(\tilde q)}$ and
$\ival{DP_{\ival c}(Q_\Kr)}$ enclose $P_c(\tilde q)$ and $DP_c(Q_\Kr)$ for $c\in\ival c$. Then for each
$c\in\ival c$ the map $P_c$ has a unique fixed point $p_c$ in $Q_\Kr$, and $p_c\in\Kr$.
\end{lemma}

(The inclusion $\Kr\subset\operatorname{int}Q_\Kr$ itself implies that $R$ and all matrices in
$\ival{DP_{\ival c}(Q_\Kr)}-I$ are invertible; see~\cite[Chapter~5]{Neumaier90}. The programs take for
$R$ an approximate inverse of the midpoint of $\ival{DP_{\ival c}(Q_\Kr)}-I$.)

Given $p_c\in\Kr_J$ for all $c\in J$, an enclosure of $p_c'=-(DP_c(p_c)-I)^{-1}\partial_cP_c(p_c)$ is
obtained by evaluating this formula in interval arithmetic on $\Kr_J\times J$, the $2\times 2$ inverse being
computed by the adjugate formula after checking that the interval determinant does not contain $0$.

\subsection{The local unstable manifold}

\begin{lemma}\label{lem:Wu}
Let $\delta,\rho,L,K,M>0$ and $B_u=[-\delta,\delta]\times[-\rho,\rho]$. Let $\ival{a_{ij}}$,
$\ival{h^i_{jk}}$, $\ival{\partial_cf^i}$ enclose the derivatives of $f$ on $B_u\times J$. Set
\begin{align*}
\Phi_1&=\ival{a_{11}}+\ival{a_{12}}[-L,L],\qquad
\Psi_1=\ival{a_{21}}+\ival{a_{22}}[-L,L],\\
\Phi_2&=\ival{h^1_{11}}+2\ival{h^1_{12}}[-L,L]+\ival{h^1_{22}}[0,L^2]+\ival{a_{12}}[-K,K],\\
\Psi_2&=\ival{h^2_{11}}+2\ival{h^2_{12}}[-L,L]+\ival{h^2_{22}}[0,L^2]+\ival{a_{22}}[-K,K].
\end{align*}
Assume:
\begin{enumerate}[label=\textup{(U\arabic*)},leftmargin=*]
\item\label{U1} $\mig\Phi_1>1$;
\item\label{U2} $L\delta\le\rho$;
\item\label{U3} $\Psi_1/\Phi_1\subset[-L,L]$;
\item\label{U4} $(\Psi_2\Phi_1-\Psi_1\Phi_2)/\Phi_1^3\subset[-K,K]$;
\item\label{U5} $\ival{\partial_cf^2}+\ival{a_{22}}[-M,M]-(\Psi_1/\Phi_1)\bigl(\ival{\partial_cf^1}+\ival{a_{12}}[-M,M]\bigr)\subset[-M,M]$;
\item\label{U6} with $m=\mig\ival{a_{11}}$ we have $1-\magn\ival{a_{12}}L/m>0$ and
$\theta:=\magn\ival{a_{21}}\,\kappa+\magn\ival{a_{22}}(L\kappa+1)<1$, where
$\kappa=\dfrac{\magn\ival{a_{12}}/m}{1-\magn\ival{a_{12}}L/m}$.
\end{enumerate}
Then for every $c\in J$ there is a function $w_c\colon[-\delta,\delta]\to[-\rho,\rho]$, real analytic in
$(u,c)$, such that:
\begin{enumerate}[label=\textup{(\alph*)}]
\item $w_c(0)=0$, $|w_c'|\le L$, $|w_c''|\le K$ and $|\partial_c w_c|\le M$;
\item $p_c+E(u,w_c(u))^T\in\Wu(p_c)$ for $|u|\le\delta$;
\item if $\ival{a^0_{ij}}$ enclose $\partial_jf^i_c(0,0)$ for $c\in J$ and
$t_0\ge\magn\bigl(\ival{a^0_{21}}/(\ival{a^0_{11}}+\ival{a^0_{12}}[-L,L]-\ival{a^0_{22}})\bigr)$, then
$|w_c'(0)|\le t_0$, and consequently $|w_c(s)|\le t_0|s|+\tfrac12Ks^2$ and $|w_c'(s)|\le t_0+K|s|$.
\end{enumerate}
\end{lemma}

The programs also check explicitly that $0\notin\ival{a_{11}}$ (which follows from~\ref{U1} anyway), so
that $m$ is the minimum of the moduli of the endpoints, and that the new graph stays in $|v|\le\rho$
(which follows from~\ref{U2}--\ref{U3}).

\begin{proof}
\emph{The class.} Let $\mathcal G_u$ be the set of functions $w\colon[-\delta,\delta]\times J\to\R$ such
that, for each $c$, $w(\cdot,c)$ is $C^1$ with $w(0,c)=0$, $|\partial_uw|\le L$ and $\partial_uw(\cdot,c)$
$K$-Lipschitz, and such that $c\mapsto w(u,c)$ is $M$-Lipschitz. By~\ref{U2}, $|w|\le L\delta\le\rho$. By
the Arzel\`a--Ascoli theorem, $\mathcal G_u$ is closed in $C^0$ with the sup norm.

\emph{Graph transform.} For $w\in\mathcal G_u$ put $\phi(u,c)=f^1_c(u,w(u,c))$ and
$\psi(u,c)=f^2_c(u,w(u,c))$; the points $(u,w(u,c))$ lie in $B_u$. Since $\partial_u\phi\in\Phi_1$,
condition~\ref{U1} makes $\phi(\cdot,c)$ monotone with $|\partial_u\phi|>1$ and $\phi(0,c)=0$, so
$\phi([-\delta,\delta],c)\supset[-\delta,\delta]$. Define $(Tw)(\cdot,c)=\psi(\cdot,c)\circ\phi(\cdot,c)^{-1}$
on $[-\delta,\delta]$; its graph is the image of the graph of $w$, restricted to $|u|\le\delta$. If $w$ is
$C^2$, then
\[
(Tw)'\circ\phi=\frac{\partial_u\psi}{\partial_u\phi},\qquad
(Tw)''\circ\phi=\frac{\partial_u^2\psi\,\partial_u\phi-\partial_u\psi\,\partial_u^2\phi}{(\partial_u\phi)^3},\qquad
\partial_c(Tw)\circ\phi=\partial_c\psi-(Tw)'\circ\phi\cdot\partial_c\phi ,
\]
with $\partial_u\psi\in\Psi_1$, $\partial_u^2\phi\in\Phi_2$, $\partial_u^2\psi\in\Psi_2$,
$\partial_c\phi\in\ival{\partial_cf^1}+\ival{a_{12}}[-M,M]$ and
$\partial_c\psi\in\ival{\partial_cf^2}+\ival{a_{22}}[-M,M]$. By~\ref{U3}--\ref{U5}, $T$ maps the analytic
elements of $\mathcal G_u$ to analytic elements of $\mathcal G_u$.

\emph{Contraction.} Let $w_1,w_2\in\mathcal G_u$ and $x=\phi_1(u_1)=\phi_2(u_2)$. The mean value theorem
applied to $f^1$ along the segment between $(u_1,w_1(u_1))$ and $(u_2,w_2(u_2))$ gives
$m|u_1-u_2|\le\magn\ival{a_{12}}(L|u_1-u_2|+\|w_1-w_2\|)$, so $|u_1-u_2|\le\kappa\|w_1-w_2\|$. Applied to
$f^2$ it gives
\[
|Tw_1(x)-Tw_2(x)|\le\magn\ival{a_{21}}|u_1-u_2|+\magn\ival{a_{22}}\bigl(L|u_1-u_2|+\|w_1-w_2\|\bigr)\le\theta\|w_1-w_2\| ,
\]
on all of $\mathcal G_u$, with $\theta<1$ by~\ref{U6}.

\emph{Fixed point.} Starting from $w\equiv0$, the iterates $T^nw$ are analytic elements of $\mathcal G_u$
forming a Cauchy sequence; their limit $w_*$ lies in $\mathcal G_u$ and satisfies $Tw_*=w_*$. For $u$ in the
domain, $f_c(\phi^{-1}(u),w_*(\phi^{-1}(u)))=(u,w_*(u))$ with $|\phi^{-1}(u)|\le|u|/\mig\Phi_1$, so the
backward orbit of every point of the graph stays on the graph and converges to $0$. Hence the graph lies
in $\Wu(0;f_c)$, which proves~(b).

\emph{Regularity.} Let $\gamma(t,c)$ be an analytic parametrisation of a local unstable manifold
$\Wu_\eta(0;f_c)$ depending analytically on $c$ (see e.g.~\cite{CFdlL03}). Near $u=0$ the graph coincides
with $\Wu_\eta(0;f_c)$, both being graphs over $u$ of points whose backward orbits stay close to $0$. For
$n$ large, the graph equals $f^n_c(\{(u,w_*(u)):|u|\le\eta'\})\cap\{|u|\le\delta\}$, so
$w_*(u,c)=\pi_2f^n_c(\gamma(t,c))$ where $t=t(u,c)$ solves $\pi_1f^n_c(\gamma(t,c))=u$. Since $\pi_1$ is
strictly monotone along the graph by~\ref{U1}, the analytic implicit function theorem shows that $w_*$ is
analytic in $(u,c)$, and its Lipschitz bounds become the derivative bounds in~(a). The saddle $p_c$, its local
unstable manifold and the parametrisation $\gamma$ continue analytically to $c$ in an open interval containing
$J$, and the non-degeneracy $\partial_t\pi_1f^n_c(\gamma(t,c))\neq0$ is an open condition that holds on the
compact set where it is used; so the same formula defines an analytic extension of $w_*$ to a neighbourhood of
$[-\delta,\delta]\times J$.

\emph{Part~(c).} The graph is tangent at $0$ to the unstable eigenvector $(1,t)$ of $A=Df_c(0,0)$, so
$t=w_c'(0)$ satisfies $t(a_{11}+a_{12}t-a_{22})=a_{21}$ with $|t|\le L$, which gives the bound. The two
consequences follow by integrating $|w_c''|\le K$.
\end{proof}

\subsection{The local stable manifold}

\begin{lemma}\label{lem:Ws}
Let $r_u,r_s,L_s,K_s,M_s>0$ and $B_s=[-r_u,r_u]\times[-r_s,r_s]$. Let $\ival{a_{ij}}$, $\ival{h^i_{jk}}$,
$\ival{\partial_cf^i}$ enclose the derivatives of $f$ on $B_s\times J$. Put
\[
\Dh=\ival{a_{11}}-[-L_s,L_s]\ival{a_{21}},\qquad
Z_1=\frac{[-L_s,L_s]\ival{a_{22}}-\ival{a_{12}}}{\Dh},\qquad \zeta=\magn Z_1,
\]
\[
\Psi_s=\ival{a_{21}}[-\zeta,\zeta]+\ival{a_{22}},\qquad
Y_i=\ival{h^i_{11}}[0,\zeta^2]+2\ival{h^i_{12}}[-\zeta,\zeta]+\ival{h^i_{22}}\quad(i=1,2).
\]
Assume:
\begin{enumerate}[label=\textup{(S\arabic*)},leftmargin=*]
\item\label{S1} $\magn\ival{a_{21}}\,r_u+\magn\ival{a_{22}}\,r_s\le r_s$;
\item\label{S2} $\mig\Dh>1$;
\item\label{S3} $m_s:=\mig\ival{a_{11}}-L_s\magn\ival{a_{21}}>0$ and
$(\magn\ival{a_{12}}+L_s\magn\ival{a_{22}})\,r_s\le m_s\,r_u$;
\item\label{S4} $Z_1\subset[-L_s,L_s]$;
\item\label{S5} $\bigl([-K_s,K_s]\Psi_s^2+[-L_s,L_s]Y_2-Y_1\bigr)/\Dh\subset[-K_s,K_s]$;
\item\label{S6} $\bigl([-L_s,L_s]\ival{\partial_cf^2}+[-M_s,M_s]-\ival{\partial_cf^1}\bigr)/\Dh\subset[-M_s,M_s]$;
\item\label{S7} $\magn\ival{a_{21}}\,L_s+\magn\ival{a_{22}}<1$.
\end{enumerate}
Then for every $c\in J$ there is $\sigma_c\colon[-r_s,r_s]\to[-r_u,r_u]$, real analytic in $(v,c)$, with
$\sigma_c(0)=0$, $|\sigma_c'|\le L_s$, $|\sigma_c''|\le K_s$ and $|\partial_c\sigma_c|\le M_s$, such that
$\{(\sigma_c(v),v):|v|\le r_s\}\subset\Ws(0;f_c)$.
\end{lemma}

\begin{proof}
Let $\mathcal G_s$ be the analogous class of functions $\sigma$ with $\sigma(0)=0$ and the Lipschitz versions
of the bounds. For $|v|\le r_s$, let $h(u)=f^1_c(u,v)-\sigma(f^2_c(u,v))$ for $|u|\le r_u$. By~\ref{S1},
$|f^2_c(u,v)|\le r_s$, so $h$ is defined; $h'=a_{11}-\sigma'a_{21}$, so $|h'|\ge m_s>0$ by~\ref{S3}; and
$|h(0)|\le(\magn\ival{a_{12}}+L_s\magn\ival{a_{22}})r_s\le m_s\,r_u$. So $h$ has a unique zero $u$ with
$|u|\le r_u$; call it $(T_s\sigma)(v)$. Clearly $T_s\sigma(0)=0$, and differentiating
$f^1(T_s\sigma(v),v)=\sigma(f^2(T_s\sigma(v),v))$ gives
\begin{align*}
(T_s\sigma)'&=\frac{\sigma'a_{22}-a_{12}}{a_{11}-\sigma'a_{21}},\qquad
\partial_c(T_s\sigma)=\frac{\sigma'\partial_cf^2+\partial_c\sigma-\partial_cf^1}{a_{11}-\sigma'a_{21}},\\
(T_s\sigma)''&=\frac{\sigma''\,(F_2')^2+\sigma'\,Y_2^*-Y_1^*}{a_{11}-\sigma'a_{21}},
\end{align*}
where $F_2'=a_{21}(T_s\sigma)'+a_{22}$ and $Y_i^*=f^i_{uu}((T_s\sigma)')^2+2f^i_{uv}(T_s\sigma)'+f^i_{vv}$, the
derivatives of $f$ being evaluated at $(T_s\sigma(v),v)\in B_s$. By~\ref{S4}--\ref{S6}, $T_s$ maps the
analytic elements of $\mathcal G_s$ into $\mathcal G_s$. For the contraction, with $u_i=T_s\sigma_i(v)$,
\[
\mig\ival{a_{11}}\,|u_1-u_2|\le L_s\magn\ival{a_{21}}|u_1-u_2|+\|\sigma_1-\sigma_2\|,
\]
so $\|T_s\sigma_1-T_s\sigma_2\|\le m_s^{-1}\|\sigma_1-\sigma_2\|$, and $m_s\ge\mig\Dh>1$ by~\ref{S2}: the interval $[-L_s,L_s]\ival{a_{21}}$ is the symmetric interval
$[-L_s\magn\ival{a_{21}},\,L_s\magn\ival{a_{21}}]$, enlarged by outward rounding, so $\mig\Dh$ is at most
$\mig\ival{a_{11}}-L_s\magn\ival{a_{21}}=m_s$ (with equality in exact arithmetic). Let
$\sigma_*$ be the fixed point, obtained as in Lemma~\ref{lem:Wu}. Then $f_c$ maps $(\sigma_*(v),v)$ to
$(\sigma_*(v_1),v_1)$ with $|v_1|\le(\magn\ival{a_{21}}L_s+\magn\ival{a_{22}})|v|$, so by~\ref{S7} forward
orbits on the graph converge to $0$. Analyticity, the analytic extension to a neighbourhood of
$[-r_s,r_s]\times J$, and the derivative bounds follow as in Lemma~\ref{lem:Wu}, using the analytic local
inverse of $f_c$.
\end{proof}

\subsection{The tangency function}

Let $N\ge1$, let $S=[s_m-\Delta,s_m+\Delta]\subset(0,\delta]$, and let $w_c$, $\sigma_c$ be given by
Lemmas~\ref{lem:Wu} and~\ref{lem:Ws}. For $s\in S$ and $c\in J$ define
\begin{equation}\label{eq:defg}
X_0(s,c)=p_c+E\bigl(s,w_c(s)\bigr)^T,\qquad \Gamma(s,c)=P_c^N(X_0(s,c)),\qquad
\begin{pmatrix}U\\V\end{pmatrix}=E^{-1}\bigl(\Gamma(s,c)-p_c\bigr),
\end{equation}
and, whenever $|V|\le r_s$, $g(s,c)=U-\sigma_c(V)$. By Lemma~\ref{lem:Wu}(b), $\Gamma(s,c)\in\Wu(p_c)$, and
$\partial_s\Gamma\neq0$ because $P_c^N$ is a local diffeomorphism. The function
$\ell_c(q)=U(q)-\sigma_c(V(q))$ is defined and real analytic on the strip
$\{p_c+E(u,v)^T:\ |v|<r_s\}$, with $\partial_u\ell_c=1$, and its zero set there is the graph
$\{(\sigma_c(v),v)\}$, an arc of $\Ws(p_c)$ by Lemma~\ref{lem:Ws}. Thus $g$ has exactly the form required in
Definition~\ref{def:tangency}; only the condition $|V|<r_s$ on $S\times J$ has to be checked. Since $w_c$ and
$\sigma_c$ extend analytically (Lemmas~\ref{lem:Wu} and~\ref{lem:Ws}), so does $g$, to a neighbourhood of
$S\times J$.

\begin{lemma}[chain rule]\label{lem:chain}
Write $D_N=DP_c^N(X_0)$, $H_N=D^2P_c^N(X_0)$ and $\partial_cP^N=\partial_cP_c^N(X_0)$: derivatives of the
$N$-th iterate of the Poincar\'e map, including the return-time corrections. Then
\begin{align*}
\partial_s g&=U_s-\sigma_c'(V)V_s,\qquad
\partial^2_sg=U_{ss}-\sigma_c''(V)V_s^2-\sigma_c'(V)V_{ss},\\
\partial_cg&=U_c-(\partial_c\sigma_c)(V)-\sigma_c'(V)V_c,
\end{align*}
where
\begin{align*}
(U_s,V_s)^T&=E^{-1}D_N\bigl(\eu+w_c'(s)\es\bigr),\\
(U_{ss},V_{ss})^T&=E^{-1}\Bigl(H_N[\eu,\eu]+2w_c'H_N[\eu,\es]+(w_c')^2H_N[\es,\es]+w_c''(s)\,D_N\es\Bigr),\\
(U_c,V_c)^T&=E^{-1}\Bigl((D_N-I)\,p_c'+\partial_cw_c(s)\,D_N\es+\partial_cP^N\Bigr).
\end{align*}
\end{lemma}

\begin{proof}
This is the chain rule, using $\partial_sX_0=\eu+w_c'\es$, $\partial_s^2X_0=w_c''\es$ and
$\partial_cX_0=p_c'+(\partial_cw_c)\es$; the term $-p_c'$ in $(U_c,V_c)$ comes from the moving origin.
\end{proof}

\subsection{The tangency criterion}

\begin{lemma}\label{lem:criterion}
Let $S=[s_-,s_+]$, $J=[c_L,c_R]$, and let $g\colon S\times J\to\R$ be continuous, twice differentiable in
$s$ and differentiable in $c$. Assume:
\begin{enumerate}[label=\textup{(T\arabic*)},leftmargin=*]
\item\label{T1} $\partial_s^2g<0$ on $S\times J$;
\item\label{T2} $\partial_sg(s_-,c)>0>\partial_sg(s_+,c)$ for all $c\in J$;
\item\label{T3} $g(s,c_L)<0$ for all $s\in S$, and $g(s_0,c_R)>0$ for some $s_0\in S$;
\item\label{T4} $\partial_cg>0$ on $S\times J$.
\end{enumerate}
Then for each $c\in J$, $\partial_sg(\cdot,c)$ has a unique zero $s_*(c)\in(s_-,s_+)$, the maximiser of
$g(\cdot,c)$; the function $G(c)=\max_Sg(\cdot,c)$ is continuous and strictly increasing with a unique
zero $\cs\in(c_L,c_R)$; at $(s_*(\cs),\cs)$ we have $g=\partial_sg=0$, $\partial_s^2g<0$, $\partial_cg>0$; and
this is the only point of $S\times J$ where $g=\partial_sg=0$. If $g$ is $C^2$, then $s_*$ and $G$ are $C^1$
with $G'(c)=\partial_cg(s_*(c),c)>0$.
\end{lemma}

\begin{proof}
By~\ref{T1} and~\ref{T2}, $\partial_sg(\cdot,c)$ is strictly decreasing and changes sign on $S$, so it has a
unique zero $s_*(c)$ in the interior, the maximiser of $g(\cdot,c)$. By compactness, $G$ is continuous.
By~\ref{T4} and the mean value theorem in $c$, $g(s,c_2)>g(s,c_1)$ for every $s$ when $c_1<c_2$, so $G$ is
strictly increasing; by~\ref{T3}, $G(c_L)<0<G(c_R)$, and the intermediate value theorem gives the unique
zero. If $g(s,c)=\partial_sg(s,c)=0$, then $s=s_*(c)$ and $G(c)=0$, so $c=\cs$. When $g$ is $C^2$, the
implicit function theorem shows that $s_*$ is $C^1$ and $G'=\partial_cg+\partial_sg\cdot s_*'=\partial_cg(s_*(c),c)$.
\end{proof}

This lemma, except for its final optional regularity statement, has been formalised in Lean~4
(Section~\ref{sec:lean}). The formal version assumes derivatives in the two-sided sense at every point of the
closed rectangle $S\times J$; this holds for our $g$, which extends analytically to a neighbourhood of
$S\times J$.

\begin{remark}[how~\ref{T3} is checked]\label{rem:taylor}
For $s\in S$, Taylor's formula and~\ref{T1} give
\[
g(s,c_L)\le g(s_m,c_L)+\partial_sg(s_m,c_L)(s-s_m)-\tfrac12\gamma(s-s_m)^2
\le g(s_m,c_L)+\frac{(\partial_sg(s_m,c_L))^2}{2\gamma},
\]
where $\gamma=\min_S|\partial_s^2g(\cdot,c_L)|$. It therefore suffices to enclose $g$ and $\partial_sg$ at the
single point $(s_m,c_L)$, and $\partial_s^2g$ on $S\times\{c_L\}$. This bound is also formalised in Lean.
\end{remark}

\section{The computer-assisted verification near $c=5.7$}\label{sec:cap}

In this section we prove the following result, and deduce the Main Theorem.

\begin{proposition}[computer-assisted]\label{prop:tangency}
For $i=1,2,3$ let $J_i=[c^{(i)}_L,c^{(i)}_R]$ be the intervals with the double-precision endpoints
given in~\eqref{eq:Ji}. Then
\[
\begin{gathered}
J_1\subset(5.6994199724,\,5.6994199925),\qquad J_2\subset(5.6995796078,\,5.6995796081),\\
J_3\subset(5.70005464309,\,5.70005464312),
\end{gathered}
\]
so that $|c-5.7|<5.81\cdot10^{-4}$, $<4.21\cdot10^{-4}$ and $<5.47\cdot10^{-5}$ on $J_1$, $J_2$, $J_3$
respectively. For each $i$:
\begin{enumerate}[label=\textup{(\roman*)}]
\item For every $c\in J_i$ the map $P_c$ has a unique fixed point $p_c$ in the box $Q_{\Kr,i}$, and $p_c$
lies in the box $Q_i$ (both boxes are given in Table~\ref{tab:near57}). It is a dissipative flip saddle:
$\lambda_u\in[-2.4040,-2.4037]$, $-1<\lambda_s<0$ and $|\lambda_u\lambda_s|<2\cdot10^{-8}$.
\item There is $c_{*,i}\in J_i$ such that $\Wu(p_{c_{*,i}})$ and $\Ws(p_{c_{*,i}})$ have a quadratic
homoclinic tangency which unfolds generically as $c$ crosses $c_{*,i}$ (Definition~\ref{def:tangency}).
The tangency lies on the branch $P^{N_i}_{c}(\Wu_\loc(p_c))$ of the unstable manifold, with
$N_1=15$, $N_2=26$, $N_3=25$, within distance $4\cdot10^{-8}$ of $p_{c_{*,i}}$.
\end{enumerate}
\end{proposition}

\subsection{Numerical exploration (non-rigorous)}
\begin{figure}[t]
\centering
\includegraphics[width=0.66\textwidth]{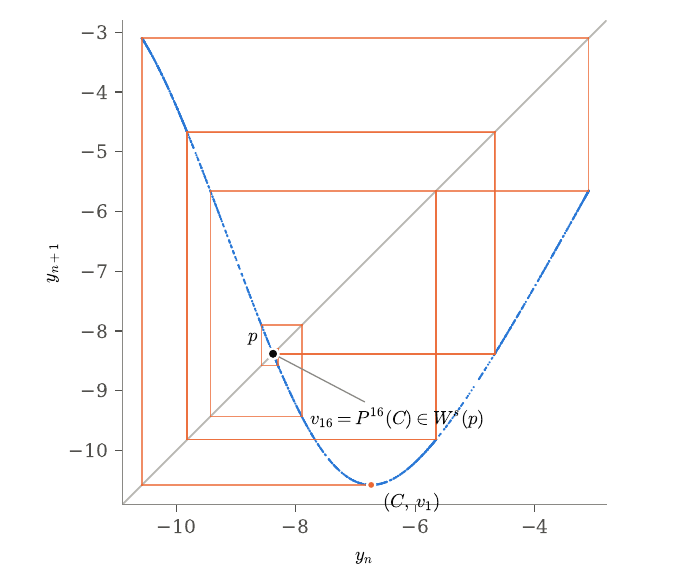}
\caption{The return map at $c=c_0\approx c_{*,3}\approx5.70005$: pairs $(y_n,y_{n+1})$ along an orbit on the
attractor (dots), and the cobweb of the critical orbit $C\mapsto v_1\mapsto\cdots\mapsto v_{16}$. At
$c_{*,3}$ the image $v_{16}$ lands on the stable manifold of the saddle $p$: this is the homoclinic
tangency of Proposition~\ref{prop:tangency} for $J_3$, used in the Main Theorem. Non-rigorous
double-precision computation.}
\label{fig:returnmap}
\end{figure}

Near $c=5.7$ the return map $P_c$ contracts areas by a factor of about $3\cdot10^{-14}$ per iterate
($\lambda_u\approx-2.40$, $\lambda_s\approx-1.3\cdot10^{-14}$). Its attractor is therefore an extremely thin
folded curve, and the dynamics is well described by the one-dimensional map $y_n\mapsto y_{n+1}$
(Figure~\ref{fig:returnmap}), unimodal on the attractor with a minimum at a critical point $C$. The fixed
point $p\approx(-8.381,\,0.0296)$ lies on its decreasing branch.

\begin{figure}[t]
\centering
\includegraphics[width=0.66\textwidth]{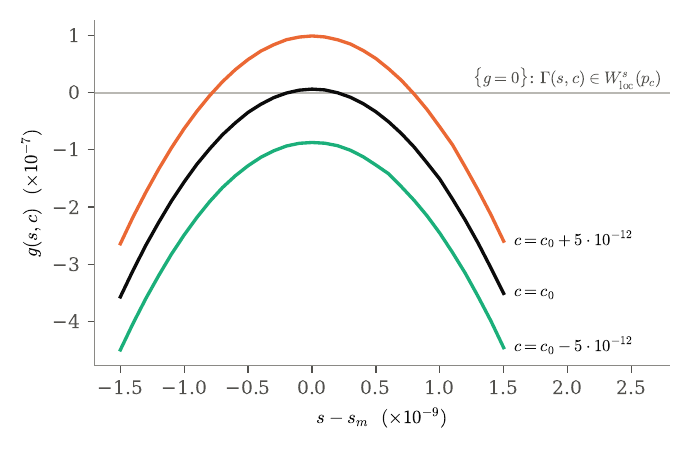}
\caption{The tangency function $g(s,c)$ for $J_3$ near its maximum, for three parameter values around the
numerical tangency parameter $c_0$. The maximum crosses $0$ (the local stable manifold) with speed
$\partial_cg\approx1.86\cdot10^4$, and $\partial_s^2g\approx-3.21\cdot10^{11}$. Non-rigorous
double-precision computation; the rigorous enclosures are in Table~\ref{tab:results}.}
\label{fig:unfolding}
\end{figure}

A homoclinic tangency of $p$ occurs when the fold of $\Wu(p)$ at $C$ is mapped onto $\Ws(p)$; in the
one-dimensional picture, when the critical point lands on $p$ after $n$ iterates. To locate such
parameters we parametrise $\Wu_\loc(p_c)$ by $s\mapsto p_c+s\eu$; $K_2=9$ or $10$ iterates bring this arc to
the fold; we then compute the $u$-coordinate of $P_c^{K_2+n}(p_c+s\eu)$ near its extremum in $s$, and find
the zeros in $c$ of the extremal value. Near $5.7$ this gives tangency parameters for many values of $n$, on
both sides of $5.7$. We selected three of them, for $n=6$ and $n=16$, close to $5.7$
(Table~\ref{tab:near57}). The matrix $E$ is the numerical eigenbasis of $DP_{c_0}(p_{c_0})$ at the
numerical tangency parameter $c_0$. Figure~\ref{fig:unfolding} shows the unfolding for $J_3$. These
approximate data only determine \emph{where} the rigorous computation looks; they play no role in the
logic of the proof.

\begin{figure}[t]
\centering
\begin{tikzpicture}[xscale=1.25,yscale=1.6]
\draw[->,gray] (-3.2,0) -- (2.4,0) node[right] {$u$};
\draw[->,gray] (0,-0.9) -- (0,1.55) node[above] {$v$};
\draw[thick,blue!70!black] (0.04,-0.85) .. controls (0.0,0) .. (-0.03,1.45);
\node[blue!70!black,anchor=west] at (0.12,1.38) {$\Ws_\loc(p_c)=\{g=0\}$};
\draw[thick] (-3.0,-0.05) .. controls (-1,0.0) and (1,0.02) .. (2.2,0.08);
\node[anchor=north] at (-2.3,-0.06) {$\Wu_\loc(p_c)$};
\fill (0,0) circle (1.6pt);
\node[anchor=north west] at (0.02,-0.02) {$p_c$};
\draw[thick,orange!85!black] (-2.9,1.0) .. controls (0.9467,0.9) and (0.9467,0.6) .. (-2.9,0.5);
\fill[orange!85!black] (-0.015,0.75) circle (1.4pt);
\node[anchor=west] at (0.08,0.75) {$q_*$};
\node[orange!85!black,anchor=south east] at (-0.35,1.02) {$\Gamma(\cdot,\cs)=P^{N}_{\cs}(\text{arc of }\Wu_\loc)$};
\draw[->,gray] (-1.25,0.3) -- (-0.45,0.3);
\node[gray,anchor=east] at (-1.3,0.3) {\small fold moves as $c$ increases};
\end{tikzpicture}
\caption{Schematic (not to scale) of the tangency in the local chart around $p_c$, with $N=25$ for $J_3$. The
fold of the arc $\Gamma(\cdot,c)$ moves with speed $\partial_cg>0$ and touches the local stable manifold
quadratically at $c=\cs$. In reality the fold is extremely thin, and the tangency point $q_*$ lies within
$4\cdot10^{-8}$ of $p_{\cs}$.}
\label{fig:schematic}
\end{figure}
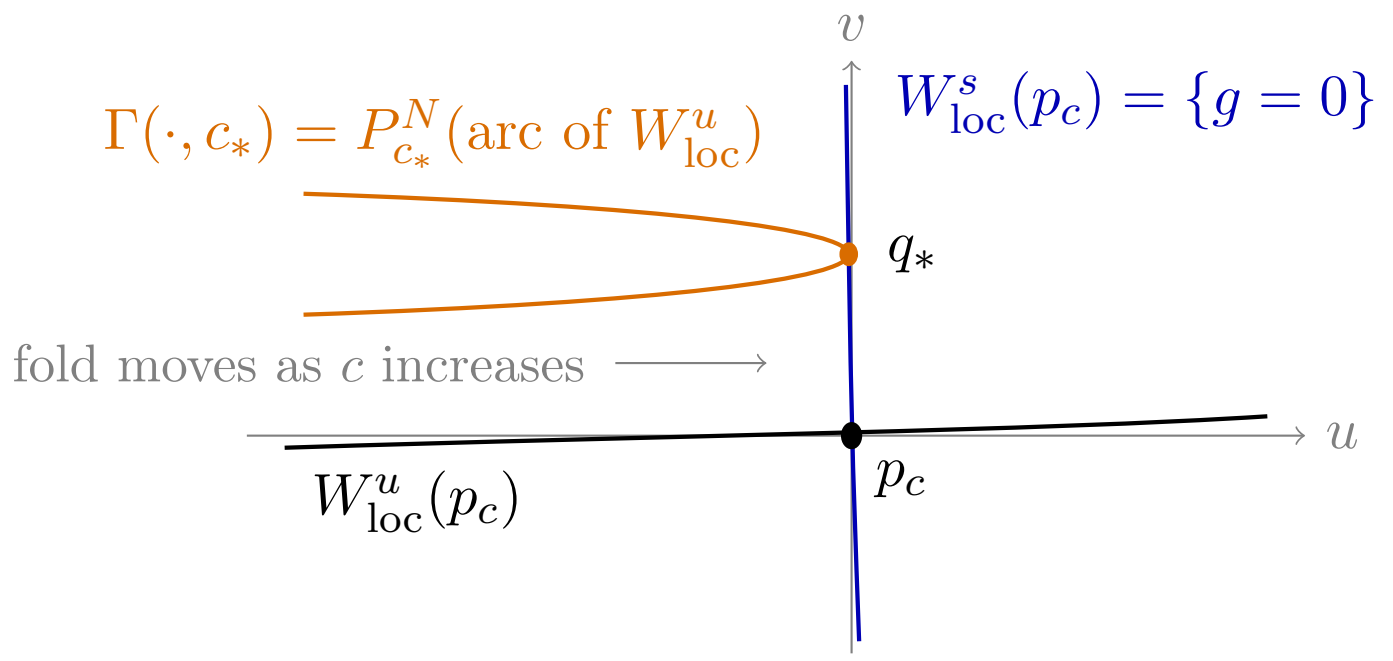

\subsection{Implementation}\label{sec:impl}
All computations in this section use the CAPD::DynSys library~\cite{CAPD21}, version 6.1.0 (GitHub
commit \texttt{2f06098}), compiled with \texttt{g++} 13.3. All the checks are implemented in one
program, \texttt{tangency\_cap.cpp}. The variants for the different intervals are generated from it by a
script that substitutes only the numerical parameters and the numerical statement being checked;
regenerating them reproduces the released files byte for byte.

\emph{Integrator.} Poincar\'e maps and their first and second derivatives are enclosed with the
$C^2$-Lohner algorithm~\cite{Lohner92,WZ11}, of Taylor order $20$, applied to the extended
four-dimensional field $\hat X$, with section $\{x=0\}$ crossed from $x<0$ to $x>0$.

\emph{Parameters.} The parameters $a=b=0.2$ enter as the intervals $\ival{2/10}$, which contain $0.2$. For
each $i$, $\cmid^{(i)}$ is the double given in Table~\ref{tab:near57}, and
$J_i=[\cmid^{(i)}-\varepsilon_c,\cmid^{(i)}+\varepsilon_c]$ with endpoints rounded to the nearest double;
exactly,
\begin{equation}\label{eq:Ji}
\begin{aligned}
J_1&=[5.6994199724407454965557917603291571140289306640625,\\
&\qquad 5.69941999244074537500637234188616275787353515625],\\
J_2&=[5.69957960782199535287873004563152790069580078125,\\
&\qquad 5.699579608021995369426804245449602603912353515625],\\
J_3&=[5.70005464309744258599721433711238205432891845703125,\\
&\qquad 5.70005464311744258765202175709418952465057373046875].
\end{aligned}
\end{equation}
The programs check that $c_L,\cmid,c_R\in J_i$.

\emph{Fixed point.} The Krawczyk test of Lemma~\ref{lem:krawczyk} is performed for $c\in J$ on the box
$\tilde q+[\pm3\cdot10^{-8}]\times[\pm10^{-9}]$, and at the single values $c_L,\cmid,c_R$ on boxes of radius
$10^{-11}$ centred at $\tilde q+\tilde p'(c-\cmid)$, with $\tilde p'$ an approximation of $p_c'$. The images
$\Kr_L,\Kr_{\rm mid},\Kr_R$ of the thin tests are checked to lie in $\Kr_J$, so all four describe the same
branch $p_c$.

\emph{Derivative boxes.} As explained in Section~\ref{sec:notation}, the derivatives used in
Lemmas~\ref{lem:Wu} and~\ref{lem:Ws} are enclosed on the interval hulls of $\Kr_J+E\cdot B_u$ and
$\Kr_J+E\cdot B_s$, with $c\in J$. The enclosures of $Df_c(0,0)$ used in Lemma~\ref{lem:Wu}(c) are taken on
$\Kr_J\times J$, and on the thin enclosures $\Kr_L\times\{c_L\}$ and $\Kr_R\times\{c_R\}$.

\emph{Second derivatives.} CAPD returns second-order Taylor coefficients, which we convert to second
derivatives by doubling the diagonal entries (Section~\ref{sec:audit}).

\emph{Initial sets.} For the composite map $P_c^N$ on $S\times J$, the initial set is an affine (Lohner
doubleton) set $\hat x+A\,r_0+r$ in the extended phase space, where $\hat x=(0,\tilde p+s_m\eu,\cmid)$ is a
point, $\tilde p$ the midpoint of the thin Krawczyk enclosure $\Kr_{\rm mid}$ of $p_{\cmid}$, the columns of
$A$ are $\eu$, $\es$ and $(\tilde p',1)$, with $\tilde p'$ the midpoint of the enclosure $\ival{p'}$ of $p_c'$,
and $r_0=(S-s_m,\ \ival{w},\ J-\cmid)$ with $\ival w$ a bound for $w_c(S)$ from Lemma~\ref{lem:Wu}(c). The
box $r$ absorbs rounding errors of $\hat x$ and the remainder
$p_c-\tilde p-\tilde p'(c-\cmid)\in(\Kr_{\rm mid}-\tilde p)+(\ival{p'}-\tilde p')(J-\cmid)$. This set contains
$\{X_0(s,c)\}$ for all $s\in S$ and $c\in J$ and keeps the strong correlation between the base point and the
parameter; the derivative enclosures returned by CAPD are valid at every point of it. For the evaluations
at the single parameters $c_L$ and $c_R$ used for~\ref{T3}, the set is centred instead at the midpoint of
$\Kr_L$ (respectively $\Kr_R$), with no parameter direction, and $r$ contains $\Kr_L$ (respectively $\Kr_R$)
minus that midpoint. We then evaluate the formulas of Lemma~\ref{lem:chain} in interval arithmetic, keeping
$D_N\eu$ and $D_N\es$ as single vectors to avoid spurious cancellation errors.

\emph{Proof parameters.} They are listed in Table~\ref{tab:near57}; in all cases $\rho=10^{-5}$ and
$r_u=r_s=10^{-6}$. The bounds $L,K,M,L_s,K_s,M_s$ are found by iterating the right-hand sides
of~\ref{U3}--\ref{U5} and~\ref{S4}--\ref{S6} with a safety factor $1.05$, and are then verified.

\emph{The statement itself.} Besides the hypotheses of Lemmas~\ref{lem:krawczyk}--\ref{lem:criterion}, the
programs check the numerical statement of Proposition~\ref{prop:tangency} as printed: that the decimal box $Q_{\Kr,i}$
lies in the Krawczyk box, that $\Kr_J\subset Q_i$, the range of $\lambda_u$ (from
$\lambda_u=\operatorname{tr}-\lambda_s$ and $|\lambda_s|\le|\det|/(\mig\operatorname{tr}-1)$), the bound on
$|\det|$, and the bound on the distance of the tangency point to $p$ (see the proof below). Decimal
constants are converted to intervals exactly enclosing them.

\subsection{Results}
Each program performs $40$ checks ($35$ for the hypotheses and $5$ for the statement); all succeed for
$J_1$, $J_2$ and $J_3$. Table~\ref{tab:results} gives the main enclosures for $J_3$, the interval closest
to $5.7$, and Table~\ref{tab:near57} summarises all three.

\begin{table}[t]
\centering
\footnotesize
\begin{tabular}{@{}lll@{}}
\toprule
Step & Quantity & Rigorous enclosure / bound\\
\midrule
Lemma~\ref{lem:krawczyk} & $p_c$, $c\in J_3$ & $[-8.38100671276265,-8.38100671273825]$\\
& & $\quad\times[0.0295898183755,0.0295898183762]$\\
& $p_c'$ & $[-1.18898577,-1.18898576]\times[-0.0044333646,-0.0044333645]$\\
saddle & $\operatorname{tr}Df_c(0)$ & $[-2.4039731272,-2.4039731256]$\\
& $\det Df_c(0)$ & $[-1.97\cdot10^{-11},\,1.98\cdot10^{-11}]$\\
Lemma~\ref{lem:Wu} & $\mig\Phi_1$ & $>2.39$\\
& $L,\ K,\ M$ & $\le9.74\cdot10^{-5},\ \le3.74\cdot10^{-4},\ \le6.35\cdot10^{-4}$\\
& $\theta$ & $\le2.07\cdot10^{-4}$\\
& $t_0$ (on $J_3$; at $c_L,c_R$) & $\le3.69\cdot10^{-12}$;\ $\le3.83\cdot10^{-14}$\\
Lemma~\ref{lem:Ws} & $L_s,\ K_s,\ M_s$ & $\le2.85\cdot10^{-5},\ \le0.714,\ \le1.33\cdot10^{-4}$\\
& $(\mig\Dh)^{-1}$ & $<0.4160$\\
\ref{T3} & $g(s_m,c_L)$ & $[-2.0093\cdot10^{-7},-1.7142\cdot10^{-7}]$\\
& $\partial_sg(s_m,c_L)$ & $[0.997,\,1.834]$\\
& $\sup_Sg(\cdot,c_L)$ & $\le-1.7142\cdot10^{-7}$\\
& $g(s_m,c_R)$ & $[1.7143\cdot10^{-7},\,2.0070\cdot10^{-7}]$\\
\ref{T2} & $\partial_sg(s_m-\Delta,c)$, $c\in J_3$ & $[4.137,\,5.107]$\\
& $\partial_sg(s_m+\Delta,c)$, $c\in J_3$ & $[-2.281,\,-1.313]$\\
\ref{T1} & $\partial_s^2g$ on $S\times J_3$ & $[-3.21014\cdot10^{11},\,-3.21010\cdot10^{11}]$\\
\ref{T4} & $\partial_cg$ on $S\times J_3$ & $[1.8564\cdot10^{4},\,1.8657\cdot10^{4}]$\\
chart & $|V|$ on $S\times J_3$ & $\le2.43\cdot10^{-9}$\ (required: $<r_s=10^{-6}$)\\
& $|U|$ on $S\times J_3$ & $\le2.21\cdot10^{-7}$\ (not required)\\
statement & $\lambda_u$; $|q_*-p|$ & $[-2.403973128,-2.403973125]$; $\le2.43\cdot10^{-9}$\\
\bottomrule
\end{tabular}
\vspace{4pt}
\caption{Rigorous enclosures for $J_3$ (the tangency closest to $c=5.7$) produced by the CAPD program, rounded
outwards. Here $f_c$ is the local map in the coordinates $(u,v)$ centred at $p_c$.}
\label{tab:results}
\end{table}

\begin{table}[t]
\centering
\scriptsize
\setlength{\tabcolsep}{2.4pt}
\begin{tabular}{@{}llll@{}}
\toprule
& $J_1$ ($n=6$) & $J_2$ ($n=16$) & $J_3$ ($n=16$)\\
\midrule
$\cmid$ & $5.6994199824407454$ & $5.6995796079219954$ & $5.7000546431074426$\\
$\varepsilon_c$, $\Delta$ & $10^{-8}$, $10^{-9}$ & $10^{-10}$, $10^{-10}$ & $10^{-11}$, $10^{-11}$\\
$K_2$, $N$, $\delta$ & $9$, $15$, $7\cdot10^{-4}$ & $10$, $26$, $6\cdot10^{-4}$ & $9$, $25$, $7\cdot10^{-4}$\\
$s_m$ & $5.9787937263417563\cdot10^{-4}$ & $2.486834987410949\cdot10^{-4}$ & $5.9741814151543799\cdot10^{-4}$\\
$\tilde q$ & $(-8.3802521043854075,$ & $(-8.3804418993899592,$ & $(-8.3810067127504517,$\\
& $\ \ 0.029592632330554456)$ & $\ \ 0.029591924532709056)$ & $\ \ 0.029589818375831582)$\\
$Q_\Kr$: $y$ & $[-8.3802521343,-8.3802520744]$ & $[-8.3804419293,-8.3804418694]$ & $[-8.3810067427,-8.3810066828]$\\
$Q_\Kr$: $z$ & $[0.029592631331,0.029592633330]$ & $[0.029591923533,0.029591925532]$ & $[0.029589817376,0.029589819375]$\\
$Q$: $y$ & $[-8.380252117,-8.380252092]$ & $[-8.38044189952,-8.38044189926]$ & $[-8.38100671277,-8.38100671273]$\\
$Q$: $z$ & $[0.0295926320,0.0295926327]$ & $[0.0295919245296,0.0295919245358]$ & $[0.0295898183755,0.0295898183762]$\\
$\operatorname{tr}Df_c(0)$ & $[-2.40374631,-2.40374475]$ & $[-2.403802784,-2.403802768]$ & $[-2.4039731272,-2.4039731256]$\\
$|\det Df_c(0)|$ & $<1.95\cdot10^{-8}$ & $<1.96\cdot10^{-10}$ & $<1.98\cdot10^{-11}$\\
$g(s_m,c_L)$ & $[-9.2878,-9.2758]\cdot10^{-8}$ & $[-2.8412,-2.8072]\cdot10^{-6}$ & $[-2.0093,-1.7142]\cdot10^{-7}$\\
$g(s_m,c_R)$ & $[9.2758,9.2878]\cdot10^{-8}$ & $[2.8069,2.8412]\cdot10^{-6}$ & $[1.7143,2.0070]\cdot10^{-7}$\\
$\partial_s^2g$ & $[-1.60020,-1.59973]\cdot10^{8}$ & $[-2.81306,-2.81237]\cdot10^{12}$ & $[-3.21014,-3.21010]\cdot10^{11}$\\
$\partial_cg$ & $[4.786,13.78]$ & $[2.380,3.269]\cdot10^{4}$ & $[1.8564,1.8657]\cdot10^{4}$\\
$|V|$ & $\le2.02\cdot10^{-9}$ & $\le3.13\cdot10^{-8}$ & $\le2.43\cdot10^{-9}$\\
\bottomrule
\end{tabular}
\vspace{4pt}
\caption{The three tangencies near $c=5.7$ (CAPD). $\tilde q$ is the (exact double) centre of the Krawczyk
box; $Q_\Kr$ is a decimal box inside the Krawczyk box (uniqueness of $p_c$), and $Q$ a decimal box containing
the Krawczyk image $\Kr_J$ (location of $p_c$). Enclosures are rounded outwards; those of $\partial_s^2g$, $\partial_cg$ and $|V|$ are over $S\times J$. For $J_1$ and $J_3$ the
eigenvector $\eu$ is replaced by $-\eu$, so that $s>0$ reaches the fold after an odd number $K_2$ of
iterates.}
\label{tab:near57}
\end{table}

\begin{proof}[Proof of Proposition~\ref{prop:tangency}]
Fix $i\in\{1,2,3\}$. Lemma~\ref{lem:krawczyk}, applied as described above, gives the fixed point $p_c$, unique
in the Krawczyk box, which contains $Q_{\Kr,i}$, and lying in $\Kr_J\subset Q_i$. Let
$\chi(\lambda)=\lambda^2-\operatorname{tr}\cdot\lambda+\det$ be the characteristic polynomial of $Df_c(0)$. The
trace and determinant enclosures give $\chi(1)>0>\chi(-1)$, hence $\lambda_u<-1<\lambda_s<1$; then
$|\lambda_u\lambda_s|=|\det|<2\cdot10^{-8}$, and $\lambda_u=\operatorname{tr}-\lambda_s$ with
$|\lambda_s|\le|\det|/(\mig\operatorname{tr}-1)$ gives the stated range of $\lambda_u$. By
Section~\ref{sec:setting}, $\det>0$, so $\lambda_s<0$. (The determinant enclosure of the CAPD program contains
$0$; it is the orientation argument that fixes the sign. The independent Kv enclosures of
Section~\ref{sec:kv} are strictly positive.)

The hypotheses~\ref{U1}--\ref{U6} and~\ref{S1}--\ref{S7} are verified, so Lemmas~\ref{lem:Wu}
and~\ref{lem:Ws} provide $w_c$ and $\sigma_c$. The function $g$ of~\eqref{eq:defg} is defined on $S\times J_i$,
since $|V|<r_s$ there. Lemma~\ref{lem:chain} and the CAPD enclosures yield~\ref{T1}--\ref{T4} (for~\ref{T3}
via Remark~\ref{rem:taylor}), and Lemma~\ref{lem:criterion} gives $c_{*,i}\in(c_L,c_R)$ and $s_*$ with
$g=\partial_sg=0$, $\partial_s^2g<0$ and $\partial_cg>0$. By the discussion preceding Lemma~\ref{lem:chain}, this
is a quadratic homoclinic tangency unfolding generically, in the sense of Definition~\ref{def:tangency}.

Finally, at the tangency $U=\sigma(V)$, so $|U|\le L_s|V|$, and the tangency point $q_*=p+U\eu+V\es$
satisfies $|q_*-p|\le(L_s\|\eu\|+\|\es\|)\sup_{S\times J}|V|$, which the programs bound by
$3.13\cdot10^{-8}<4\cdot10^{-8}$ (Figure~\ref{fig:schematic}). It differs from $p$, because $q_*=P^{N_i}(X_0)$
with $X_0\neq p$ and $P$ is injective.
\end{proof}

\begin{remark}[uniqueness]
Each $c_{*,i}$ is unique only for the tracked branch $\Gamma(S,c)$. In any generic unfolding of a dissipative
homoclinic tangency, further tangencies of $\Wu(p_c)$ and $\Ws(p_c)$, created by longer iterates, accumulate
at $c_{*,i}$~\cite{PT93}.
\end{remark}

\begin{remark}[sanity checks]\label{rem:sanity}
The double-precision estimates of $\partial_s^2g$ and $\partial_cg$ quoted in Figure~\ref{fig:unfolding} lie in
the rigorous enclosures. As negative tests, we shifted each of $J_0,\dots,J_3$ by $\pm3\varepsilon_c$, so that it
no longer contains the tangency (the checks of the statement are then skipped): in all eight runs exactly
one of the remaining $35$ checks fails, namely~\ref{T3}, at $c_L$ for the positive shift and at $c_R$ for the
negative one, as predicted by $\partial_cg>0$. The derivatives of $P^N_c$ were also compared with finite
differences (Section~\ref{sec:audit}).
\end{remark}

\begin{proof}[Proof of the Main Theorem]
Let $\cs=c_{*,3}$ be the parameter given by Proposition~\ref{prop:tangency} for $i=3$. Then
$|\cs-5.7|<5.47\cdot10^{-5}$, the saddle $p_{\cs}$ satisfies $|\det DP_{\cs}(p_{\cs})|<2\cdot10^{-8}<1$, its
multipliers have no resonance of order below $1\,229\,030$ (Remark~\ref{rem:MVhyp}), and its stable and unstable
manifolds have a quadratic homoclinic tangency unfolding generically at $\cs$. The first
statement therefore follows from Proposition~\ref{prop:srb}, and the second by taking
$\varepsilon=5\cdot10^{-6}$, since $(\cs-5\cdot10^{-6},\cs+5\cdot10^{-6})\subset(5.7-6\cdot10^{-5},5.7+6\cdot10^{-5})$.
The same argument applies to $c_{*,1}$ and to the band-merging parameter of Proposition~\ref{prop:band}, for
which the non-resonance condition is also checked by the Kv programs. For $c_{*,2}$, which was verified with CAPD
only, the CAPD enclosures $|\det DP_c(p_c)|<1.96\cdot10^{-10}$ and $|\lambda_u|>2.4038$ give $\rho>26.4$, which
excludes resonances of order at most $28$ only.
\end{proof}

\section{The band-merging tangency}\label{sec:band}

\begin{proposition}[computer-assisted]\label{prop:band}
Let $J_0=[c^{(0)}_L,c^{(0)}_R]$, with the double-precision endpoints given in~\eqref{eq:J0}, so that
$J_0\subset(4.506989336,\,4.506989357)$. Then the conclusions (i)--(ii) of
Proposition~\ref{prop:tangency} hold on $J_0$, with the boxes $Q_{\Kr,0}$ and $Q_0$ given below,
$\lambda_u\in[-1.92148,-1.92146]$, $|\lambda_u\lambda_s|<5\cdot10^{-8}$, $N_0=15$, and the tangency point within
$2\cdot10^{-8}$ of $p_{c_{*,0}}$. Consequently, the conclusion of Proposition~\ref{prop:srb} holds near
$c_{*,0}$.
\end{proposition}

The same program, with the parameters of Table~\ref{tab:band}, proves Proposition~\ref{prop:band}. Here $K_2=12$ and
$N=15$: the fold of $\Wu(p)$ is mapped onto $\Ws(p)$ after $3$ further iterates. Numerically, this is the
parameter at which the two bands of the chaotic attractor of the period-doubling cascade merge. The interval
is
\begin{equation}\label{eq:J0}
\begin{aligned}
J_0=[&4.50698933631197018456759906257502734661102294921875,\\
&4.50698935631197006301817964413203299045562744140625],
\end{aligned}
\end{equation}
the doubles nearest to $\cmid\mp10^{-8}$, with $\cmid$ the double nearest to $4.50698934631197$. The Krawczyk
test gives a unique fixed point in
\[
Q_{\Kr,0}=[-6.9386987611,\,-6.9386987012]\times[0.0360634877,\,0.0360634896],
\]
lying in $Q_0=[-6.9386987439,\,-6.9386987185]\times[0.0360634882,\,0.0360634892]$. All $40$ checks succeed,
and the shifted intervals $J_0\pm3\cdot10^{-8}$ fail at~\ref{T3} only, as expected. The proof of
Proposition~\ref{prop:band} is then identical to that of Proposition~\ref{prop:tangency}; here the distance of the tangency point
to $p$ is bounded by $9.92\cdot10^{-9}<2\cdot10^{-8}$.

\begin{table}[t]
\centering
\footnotesize
\begin{tabular}{@{}lll@{}}
\toprule
Step & Quantity & Rigorous enclosure / bound\\
\midrule
parameters & $s_m$, $\Delta$, $\delta$ & $5.3682155586\cdot10^{-4}$, $10^{-9}$, $6\cdot10^{-4}$\\
& $E$ & $\begin{pmatrix}0.9999996961880596 & -0.8077389541129416\\ 7.795022697815469\cdot10^{-4} & -0.5895403141503822\end{pmatrix}$\\
saddle & $\operatorname{tr}Df_c(0)$, $\det Df_c(0)$ & $[-1.92147167,-1.92147014]$, $[-4.52,4.53]\cdot10^{-8}$\\
\ref{T3} & $g(s_m,c_L)$, $g(s_m,c_R)$ & $[-1.7531,-1.6882]\cdot10^{-8}$, $[1.6881,1.7529]\cdot10^{-8}$\\
\ref{T2} & $\partial_sg(s_m\mp\Delta,c)$ & $[0.01908,0.02497]$, $[-0.02557,-0.01967]$\\
\ref{T1} & $\partial_s^2g$ on $S\times J_0$ & $[-2.23388,-2.23090]\cdot10^{7}$\\
\ref{T4} & $\partial_cg$ on $S\times J_0$ & $[1.4164,2.0247]$\\
chart & $|V|$ on $S\times J_0$ & $\le9.92\cdot10^{-9}$\\
\bottomrule
\end{tabular}
\vspace{4pt}
\caption{The band-merging tangency of Proposition~\ref{prop:band} (CAPD), rounded outwards.}
\label{tab:band}
\end{table}

\section{Independent verification, formal proof and code}\label{sec:indep}

\subsection{A second implementation with the Kv library}\label{sec:kv}
We re-verified all hypotheses used in the proofs for $J_1$, $J_3$ and $J_0$ (the Krawczyk test,
\ref{U1}--\ref{U6}, \ref{S1}--\ref{S7}, \ref{T1}--\ref{T4} and the chart condition), together with the
numerical statements (exact endpoints of $J$, the boxes $Q_\Kr$ and $Q$, the range of $\lambda_u$, the bound
on $|\det|$, and the absence of resonances of order at most $10^3$ used in Remark~\ref{rem:MVhyp}), with a separately written program built on Kv~\cite{Kv}, a C++ library for verified numerical
computation by M.~Kashiwagi (version 0.4.62), instead of CAPD. Three features make this implementation
differ in substance, not only in code.
\begin{itemize}[leftmargin=*]
\item \emph{Poincar\'e map as a fixed-time map.} Write $x=r\cos\theta$, $y=r\sin\theta$, and use
$t=\theta+\pi/2$ as independent variable. Then
\[
\frac{dr}{dt}=\frac{r(arS^2-Cz)}{D},\qquad\frac{dz}{dt}=\frac{r(b+z(rC-c))}{D},
\]
where $C=\sin t$, $S=-\cos t$ and $D=r(1+aCS)+Sz=r\dot\theta$. The section $\Sigma$ is $\{t\equiv0\bmod2\pi\}$,
so as long as $\dot\theta>0$, $P_c$ is exactly the time-$2\pi$ map in the coordinates $(y,z)=(-r,z)$, and $P_c^N$
the time-$2\pi N$ map: no return-time correction is needed. The condition $D>0$ is enforced rigorously, since
Kv's series division fails if the range of $D$ contains $0$ on any Taylor step. (A non-rigorous exploration,
not used in the proof, gives $\dot\theta\ge0.88$ on the relevant orbits.)
\item \emph{Jets.} First and second derivatives are obtained by integrating the first and second variational
equations, written with forward-mode jet arithmetic inside the vector field, seeded along the directions
actually needed (for example $\eu+w'\es$ with $w'$ an interval). In the $c$-slot,
$D_N\tilde p'+\partial_cP^N$ is computed as a single directional derivative.
\item \emph{Sets.} All sets are Kv affine forms, propagated with the mean-value form of Kv's affine ODE solver,
with noise symbols for $s$, $c$, the $w$-range and the remainder of $p_c$.
\end{itemize}
The one-step enclosures are Kv's own; the program adds a step-size strategy and exact power-of-two rescalings
of the jet components, which do not affect rigour. It uses the same proof parameters as the CAPD programs, and
its variants for $J_1$ and $J_3$ are generated from the program for $J_0$ by substituting only the numerical
parameters and the numerical statement.

All $38$ checks ($32$ for the hypotheses and $6$ for the statement) succeed in each case, with runtimes on one core of $9$, $13$ and $7$ minutes for $J_1$,
$J_3$ and $J_0$. The Kv determinant enclosures are strictly positive, confirming $\lambda_s\in(-1,0)$
independently of Section~\ref{sec:setting}, and the Kv enclosures of $\det$ and $\partial_cg$ are much tighter
than CAPD's; the others are comparable (Table~\ref{tab:kv}). The negative tests ($J_0$ and $J_3$ shifted by
$\pm3\varepsilon_c$) fail exactly at~\ref{T3}.

The Kv program was written independently of the CAPD program at the level of numerical method, library and
code; however, it was produced with the same AI assistance (see the statement at the end) and has not yet been
audited by an independent human expert. Its value is to rule out errors specific to one library, one
integration scheme, or one implementation of the formula layer.

\begin{table}[t]
\centering
\footnotesize
\begin{tabular}{@{}llll@{}}
\toprule
& Quantity & CAPD & Kv\\
\midrule
$J_3$ & $\det Df_c(0)$ & $[-1.97,1.98]\cdot10^{-11}$ & $[3.08113,3.08114]\cdot10^{-14}$\\
& $g(s_m,c_L)$ & $[-2.0093,-1.7142]\cdot10^{-7}$ & $[-2.0885,-1.6210]\cdot10^{-7}$\\
& $g(s_m,c_R)$ & $[1.7143,2.0070]\cdot10^{-7}$ & $[1.8007,1.9247]\cdot10^{-7}$\\
& $\partial_s^2g$ on $S\times J_3$ & $[-3.21014,-3.21010]\cdot10^{11}$ & $[-3.21015,-3.21009]\cdot10^{11}$\\
& $\partial_cg$ on $S\times J_3$ & $[1.8564,1.8657]\cdot10^{4}$ & $[1.861053,1.861065]\cdot10^{4}$\\
\midrule
$J_1$ & $\det Df_c(0)$ & $[-1.95,1.95]\cdot10^{-8}$ & $[3.0897,3.0957]\cdot10^{-14}$\\
& $g(s_m,c_L)$ & $[-9.2878,-9.2758]\cdot10^{-8}$ & $[-9.28208,-9.28140]\cdot10^{-8}$\\
& $\partial_s^2g$ on $S\times J_1$ & $[-1.60020,-1.59973]\cdot10^{8}$ & $[-1.599966,-1.599959]\cdot10^{8}$\\
& $\partial_cg$ on $S\times J_1$ & $[4.786,13.78]$ & $[9.28165,9.28190]$\\
\midrule
$J_0$ & $\det Df_c(0)$ & $[-4.52,4.53]\cdot10^{-8}$ & $[3.43502,3.43503]\cdot10^{-11}$\\
& $g(s_m,c_L)$ & $[-1.7531,-1.6882]\cdot10^{-8}$ & $[-1.72066,-1.72056]\cdot10^{-8}$\\
& $\partial_s^2g$ on $S\times J_0$ & $[-2.23388,-2.23090]\cdot10^{7}$ & $[-2.2323946,-2.2323828]\cdot10^{7}$\\
& $\partial_cg$ on $S\times J_0$ & $[1.4164,2.0247]$ & $[1.720502,1.720572]$\\
\bottomrule
\end{tabular}
\vspace{4pt}
\caption{Independent verification: enclosures from the CAPD and the Kv programs, both rounded outwards.}
\label{tab:kv}
\end{table}

\subsection{Formalisation of the tangency criterion in Lean}\label{sec:lean}
Lemma~\ref{lem:criterion} and the bound of Remark~\ref{rem:taylor} have been formalised in Lean~4 (v4.35.0)
with Mathlib~\cite{Lean4,mathlib} (commit \texttt{aa35a48}). The proofs contain no \texttt{sorry} and no
additional axioms: \texttt{\#print axioms} reports only \texttt{propext}, \texttt{Classical.choice} and
\texttt{Quot.sound}. The main statement, slightly abridged, reads:
{\footnotesize
\begin{verbatim}
theorem tangency_criterion
  (hg  : ∀ s ∈ Icc sL sR, ∀ c ∈ Icc cL cR, HasDerivAt (fun s => g s c) (gs s c) s)
  (hgs : ∀ s ∈ Icc sL sR, ∀ c ∈ Icc cL cR, HasDerivAt (fun s => gs s c) (gss s c) s)
  (hgc : ∀ s ∈ Icc sL sR, ∀ c ∈ Icc cL cR, HasDerivAt (fun c => g s c) (gc s c) c)
  (hcont : ContinuousOn (fun p : ℝ × ℝ => g p.1 p.2) (Icc sL sR ×ˢ Icc cL cR))
  (hs : sL < sR) (hc : cL < cR)
  (hA : ∀ s ∈ Icc sL sR, ∀ c ∈ Icc cL cR, gss s c < 0)
  (hB : ∀ c ∈ Icc cL cR, 0 < gs sL c ∧ gs sR c < 0)
  (hC : (∀ s ∈ Icc sL sR, g s cL < 0) ∧ ∃ s₀ ∈ Icc sL sR, 0 < g s₀ cR)
  (hD : ∀ s ∈ Icc sL sR, ∀ c ∈ Icc cL cR, 0 < gc s c) :
  ∃ cstar ∈ Ioo cL cR, ∃ sstar ∈ Ioo sL sR,
    g sstar cstar = 0 ∧ gs sstar cstar = 0 ∧ gss sstar cstar < 0 ∧ 0 < gc sstar cstar ∧
    (∀ c ∈ Icc cL cR, ∃! s, s ∈ Icc sL sR ∧ gs s c = 0) ∧ ...
\end{verbatim}
}
\noindent
The omitted conclusions state that $\cs$ is the unique zero of $G$ on $[c_L,c_R]$, and that $(s_*,\cs)$ is the only
point of the rectangle where $g=\partial_sg=0$. The formal proof follows the monotonicity argument of
Lemma~\ref{lem:criterion}: it uses the mean value theorem, compactness and the intermediate value theorem, but
not the implicit function theorem.

What is formalised is the deduction from~\ref{T1}--\ref{T4} to the tangency; the graph-transform lemmas, the
chain rule, the theory of Poincar\'e maps and the interval computations are not. A formal verification of the
numerical part would require a verified validated ODE solver with second-order variational equations,
comparable to Immler's verification of Tucker's computation in Isabelle/HOL~\cite{Immler18}; this does not
currently exist in Lean.

\subsection{Code audit and release}\label{sec:audit}
The points of the implementation where an error would be most consequential were checked as follows.
\begin{itemize}[leftmargin=*]
\item \emph{Derivatives of the Poincar\'e map, not of the flow.} In the CAPD programs, $P$, $DP$ and $D^2P$ are
computed by CAPD's Poincar\'e-map routines (\texttt{computeDP}), which return the derivatives of the
Poincar\'e map itself, including the first- and second-order corrections for the varying return time. Since
$c$ is a state variable of $\hat X$, the same correction applies to $\partial_c$. For $P^N$ we integrate up to
the $N$-th crossing and apply the same routine, which yields $D(P^N)$ and $D^2(P^N)$. A separate
(non-rigorous) test program compares all first and second derivatives in the directions $\eu$, $\es$ and $c$,
including those with respect to $c$ and the Taylor-coefficient convention, with central finite differences of
$P^N_c$ at the point $X_0=\tilde p+s_m\eu$ used in the proof and at a second point $\tilde p+\frac12s_m\eu$ off the
fold. For $J_0$ ($N=15$), at both points, all entries larger than $0.1$ in absolute value agree to within
$9\cdot10^{-4}$ relative, most of them to $10^{-5}$--$10^{-8}$. For $J_3$ ($N=25$), at $X_0$ the first derivatives in
the directions $\eu$ and $c$ and all second derivatives except those in the direction $(\es,\es)$ agree to within
$6\cdot10^{-3}$, those of size at least $10^7$ (up to $4.5\cdot10^{11}$) to within $8\cdot10^{-4}$; at the second point the
second derivatives (of size up to $5\cdot10^{15}$) vary on scales below the smallest usable step, and finite
differences are inconclusive there. Entries that are small at the fold (the first derivatives along the
strongly contracted direction $\es$, and the second derivatives in the direction $(\es,\es)$) cannot be resolved
by double-precision finite differences. In the Kv programs no
correction is needed, since $P_c$ is a fixed-angle map (Section~\ref{sec:kv}); the overlap of the two sets of
enclosures (Table~\ref{tab:kv}) is therefore a further cross-check of the return-time corrections, including
the $c$-direction and the $N$-fold composition.
\item \emph{$E^{-1}$.} In both implementations, $E^{-1}$ is enclosed rigorously by the adjugate formula in
interval arithmetic, after checking that the interval determinant does not contain $0$.
\item \emph{Decisions.} Every check is decided by an interval comparison: an inequality between an endpoint and
a constant, or an inclusion of intervals. Quantities such as $\mig$ and $\magn$ are computed from endpoints
exactly. Midpoints are used only to choose the Krawczyk preconditioner and the centres of affine sets, whose
remainders are included in the sets, so they do not affect rigour.
\item \emph{Enclosure boxes.} Derivatives are enclosed on $\Kr_J+E\cdot B$, as in Section~\ref{sec:notation}.
\end{itemize}
The complete code, with build scripts, the full output logs, the negative tests and a script that rebuilds
and reruns everything, is available at \url{https://github.com/cabytheo-cloud/Rossler_SRB}.

\subsection{Trust base}
Propositions~\ref{prop:tangency} and~\ref{prop:band} rest on: the mathematical arguments of Section~\ref{sec:reduction}, of
which Lemma~\ref{lem:criterion} is machine-checked; the correctness of the CAPD library, the compiler and the
floating-point rounding modes, with, for $J_1$, $J_3$ and $J_0$, a second and independent implementation on
the Kv library; and standard analytic facts (analytic dependence of local invariant manifolds, the implicit
function theorem). Proposition~\ref{prop:srb}, and hence the Main Theorem, additionally rests on the cited theorems of Mora--Viana,
Benedicks--Young and Benedicks--Viana (or Wang--Young), whose non-resonance hypothesis~\ref{MV3} is
verified by the Kv programs (Remark~\ref{rem:MVhyp}).

\section{Discussion}\label{sec:discussion}

\subsection*{The tangency parameter $\cs$}
The Main Theorem concerns parameters \emph{near} $\cs$. Whether $P_{\cs}$ itself has an SRB measure, which would give one
for the flow at $c=\cs$ by Lemma~\ref{lem:suspension}, is open. We explain what is known and why our methods do not
decide it.

The Main Theorem gives no SRB measure at $\cs$. Its measures are carried by the small attractors of
Theorem~\ref{thm:MV}, which live in the renormalised return maps $f_\mu^{\,n+N}|U_n$ for $\mu$ in parameter sets
that accumulate at $\mu=0$ without containing it: at the tangency parameter these attractors have not yet been
created. An SRB measure at $\cs$ would therefore have to be of another kind, for instance carried by the
homoclinic set of $p_{\cs}$.

There is a one-dimensional heuristic in favour of an SRB measure. In the almost one-dimensional return map of
Figure~\ref{fig:returnmap}, the tangency at $\cs$ corresponds to the critical point landing after $16$ iterates
on the repelling fixed point. For an S-unimodal interval map this is a Misiurewicz condition, and such maps have
an absolutely continuous invariant measure with a positive Lyapunov exponent~\cite{Mis81}, the one-dimensional
analogue of an SRB measure. The analogy does not survive in dimension two. The map $P_{\cs}$ is a diffeomorphism
whose folds are infinitely many, one on each turn of $\Wu(p_{\cs})$ (the ``critical points'' of Benedicks and
Carleson~\cite{BC91}); the tangency controls the forward orbit of only one of them, which converges to
$p_{\cs}$, and says nothing about the others. In the theories of Mora--Viana~\cite{MV93} and Wang--Young~\cite{WY01}, Misiurewicz maps
appear only as singular limits, and they yield positive-measure sets of parameters near, not at, the
corresponding parameter. To our knowledge, no explicit parameter of a H\'enon-like family with non-zero
Jacobian is known to have an SRB measure.

In the other direction, there is a mechanism that rules out SRB measures on the homoclinic set of $p_c$.
By~\cite[Proposition~2.1]{New04}, if the homoclinic set of $p_c$ carries an SRB measure $\nu$, then the support
of $\nu$ contains $\Wu(p_c)$; hence, if $\Wu(p_c)$ meets the basin of a periodic sink, no such measure exists.
Unfoldings of tangencies of dissipative saddles create sinks of this kind~\cite[Lemma~2.2]{New04}, and in the
H\'enon family Kiriki, Li and Soma~\cite{KLS10} used this to obtain an open dense set of parameters, in an
interval, at which the homoclinic set of the saddle carries no SRB measure, coexisting with a set of positive
measure at which a small H\'enon-like attractor does carry one. We expect the same picture near $\cs$, with
parameters of both kinds accumulating at $\cs$; we have not checked the hypotheses of~\cite{KLS10} for our
family. In any case, a tangency is not by itself an obstruction: all known SRB measures of H\'enon-like maps
are carried by homoclinic sets containing tangencies~\cite{New04}.

Finally, the question is not accessible to the computations of this paper. We know $\cs$ only through the
interval $J_3$, of width less than $3\cdot10^{-11}$, and every property we certify holds uniformly on such an
interval. A hyperbolic sink at $\cs$ would persist on a parameter interval around $\cs$, and could in principle
be certified after enclosing $\cs$ precisely enough; but the periodic windows created by the unfolding shrink
geometrically as they approach $\cs$ without containing it, so they give no reason to expect one. Neither the
absence of sinks nor the existence of an SRB measure is an open condition, so neither can be certified by
refining $J_3$. Simulations at a floating-point approximation of $\cs$
are equally uninformative. Deciding the question at $\cs$ seems to require the same control of all critical
orbits as at $c=5.7$, with the difference that one of them is already controlled.

\subsection*{The classical parameter}
The Main Theorem shows that SRB measures do occur in the R\"ossler family, in the
same measure-theoretic sense as for the H\'enon family, arbitrarily close to tangencies at distance less than
$5.5\cdot10^{-5}$ from the classical value. It says nothing about a \emph{specific} parameter. The parameter
sets produced by Mora--Viana and Benedicks--Carleson-type exclusion arguments are Cantor-like, and certifying
that one given parameter belongs to such a set would require controlling the critical orbit for all times;
this has been done in one dimension only for special (for example Misiurewicz-type) parameters.

At $c=5.7$ two scenarios are compatible with the rigorous results available: a non-uniformly hyperbolic
attractor with an SRB measure, or a periodic sink of very long period with a tiny basin, embedded in a
quasi-attractor. Galias's periodic window at distance less than $2\cdot10^{-22}$ from the classical parameters
(in the direction $a=b$)~\cite{Gal26} shows how closely the two interleave; the extremely small widths of the
windows he found are, however, strong numerical evidence against the second scenario at the classical
parameters themselves. Our results, which produce SRB parameters at distance of order $10^{-5}$ in the
direction of $c$, are consistent with this picture: by Newhouse's phenomenon (see~\cite{New04}), the unfoldings of our tangencies
also produce periodic sinks for parameters accumulating at each $c_{*,i}$.

In non-rigorous simulations, the largest Lyapunov exponent at $c=5.7$ is about $0.071$, and the sum of the
exponents about $-5.32$. For a compactly supported invariant probability $\mu$, the identities
$\int\dot x\,d\mu=\int\dot y\,d\mu=0$ give $\int x\,d\mu=a\int z\,d\mu$, so the $\mu$-average of the sum of the
Lyapunov exponents is $\int\operatorname{div}X_c\,d\mu=a-c+a\int z\,d\mu$. This makes the strong dissipativity
quantitative, but does not exclude sinks.

\begin{remark}[observable measures]\label{rem:observable}
At $(a,b,c)=(0.2,0.2,5.7)$, the box $W=[-10.7,-2.7]\times[0.028,0.034]$ satisfies $P(W)\subset\operatorname{int}W$
(this is~\cite[Theorem~2]{CAPD21}; we re-verified it with CAPD, enclosing the decimal endpoints by exact
rational interval arithmetic). By the general theorem of Catsigeras and Enrich~\cite[Theorem~1]{CE11}, every
continuous map of a compact set of positive Lebesgue measure has \emph{SRB-like} (observable) measures; hence
$P|_W$ has such measures, supported on $\bigcap_nP^n(W)$. This weak notion does not distinguish between the
two scenarios above.
\end{remark}

\subsection*{Shilnikov homoclinic orbits}
Ding and Zheng~\cite{DZ23} derive sufficient conditions for a homoclinic orbit of Shilnikov type to a
saddle-focus, using series expansions of the solutions on its stable and unstable manifolds, and apply them to
the R\"ossler system at one parameter set. A Shilnikov loop yields horseshoes nearby (and homoclinic tangencies
in its unfoldings), but by itself it does not give SRB measures: horseshoes have zero Lebesgue measure and
support no physical measure. To reach SRB measures one needs, as here, a generically unfolding tangency of a
dissipative saddle together with the Mora--Viana and Benedicks--Viana theory. A computer-assisted proof of a
Shilnikov loop, as in~\cite{CW17}, would be a natural source of further tangencies, whose genericity for the
R\"ossler family would still have to be checked.

\needspace{5\baselineskip}
\subsection*{Open problems}
\begin{enumerate}[leftmargin=*]
\item \emph{More tangencies.} The method applies, with only numerical changes, to further tangencies of $p$ from
longer iterates, or to tangencies of periodic saddles of higher period; each successful run gives a further
positive-measure set of SRB parameters. Whether tangency parameters accumulate at $c=5.7$ is not known, and
could not be established by finitely many such runs.
\item \emph{The tangency parameters.} Does $P_{c_{*,i}}$ carry an SRB measure, for instance on the homoclinic set
of $p_{c_{*,i}}$? Alternatively, does $\Wu(p_{c_{*,i}})$ meet the basin of a sink?
\item \emph{A fixed parameter.} Proving the existence of an SRB measure at $c=5.7$ seems to require a finitely
checkable condition implying non-uniform hyperbolicity for a single strongly dissipative folding map, for
instance a two-dimensional analogue of computable Collet--Eckmann criteria. No such criterion is currently
known. Galias's results~\cite{Gal26} suggest that the classical attractor is not periodic, which is consistent
with, but does not imply, the existence of an SRB measure.
\end{enumerate}

\appendix
\section{Reproducibility}
All programs, build scripts and output logs are available at \url{https://github.com/cabytheo-cloud/Rossler_SRB}:
\begin{itemize}[leftmargin=*]
\item \texttt{tangency\_cap.cpp}: the CAPD program for Proposition~\ref{prop:band};
\item \texttt{tangency\_57\_n6.cpp} ($J_1$), \texttt{tangency\_57\_n16b.cpp} ($J_2$),
\texttt{tangency\_57\_n16a.cpp} ($J_3$): the CAPD programs for Proposition~\ref{prop:tangency}
(and the Main Theorem), generated from \texttt{tangency\_cap.cpp} by \texttt{make\_variant.py};
\item \texttt{tangency\_kv.cpp} ($J_0$), \texttt{tangency\_kv\_n6.cpp} ($J_1$),
\texttt{tangency\_kv\_n16a.cpp} ($J_3$): the Kv programs, with their headers and build script;
\item \texttt{test\_fd.cpp}: the finite-difference comparison of Section~\ref{sec:audit};
\item \texttt{trap2.cpp}: the trapping-region check of Remark~\ref{rem:observable};
\item \texttt{Criterion.lean}: the Lean formalisation.
\end{itemize}
Build and run, for example, with
\begin{verbatim}
g++ -O2 tangency_57_n16a.cpp $(capd-config --cflags --libs) -o tangency_57_n16a
./tangency_57_n16a
\end{verbatim}
Optional arguments set $\varepsilon_c$ and $\Delta$, and shift $\cmid$ for the negative tests. Every check prints
\texttt{true}, and the program ends with \texttt{ALL CHECKS PASSED: true}. The Lean file is checked with
\texttt{lake env lean} inside a Mathlib checkout at the commit given above.

\section*{AI use acknowledgement}
The author made extensive use of an AI model, Claude Opus 5.5 extra (Anthropic), in this work, including for
literature review, program writing, computations, Lean formalisation and manuscript preparation, under the
instructions of the author. The latter takes full responsibility for the content of this paper.

\end{document}